\documentclass[10pt]{amsart}
\usepackage{graphicx}

\usepackage{latexsym,amsmath}
\usepackage{amsmath,amsthm}
\usepackage{amsfonts}
\usepackage[psamsfonts]{amssymb}

\theoremstyle{plain} 
\newtheorem{theorem}{Theorem}[section]
\newtheorem{corollary}[theorem]{Corollary}

\newtheorem{lemma}[theorem]{Lemma}  

\theoremstyle{definition} 

\theoremstyle{definition} 
\newtheorem{example}{Example}
\theoremstyle{remark} 

\theoremstyle{remark} 
\newtheorem{remark}[theorem]{Remark}
\newtheorem*{remark*}{Remark}
\numberwithin{equation}{section}

\newcommand{\FF}[1]{\mathcal{F}^{#1}}
\renewcommand{\H}[1]{\mathcal{H}_+^{#1}}
\renewcommand{\r}{\mathsf{r}}

\newcommand{\lc}{\mathsf{L\!C}}

\renewcommand{\le}{\leqslant}
\renewcommand{\ge}{\geqslant}

\newcommand{\lp}{\left(}
\newcommand{\rp}{\right)}

\newcommand{\al}{\alpha}

\newcommand{\si}{\sigma}
\newcommand{\la}{\lambda}
\newcommand{\de}{\delta}
\newcommand{\vpi}{\varphi}
\renewcommand{\Psi}{\overline{\Phi}}

\renewcommand{\P}{\mathsf{P}} 

\newcommand{\E}{\mathsf{E}}
\newcommand{\Var}{\mathsf{Var}}

\newcommand{\R}{\mathbb{R}}

\newcommand{\LL}{\mathcal{L}}
\newcommand{\G}{\mathcal{G}}
\newcommand{\EE}{\mathcal{E}}
\newcommand{\I}{\mathcal{I}}
\newcommand{\J}{\mathcal{J}}

\renewcommand{\d}{\mathrm{d}}

\newcommand{\vp}{\varepsilon}

\newcommand{\iii}{\operatorname{\mathbf{I}}}
\newcommand{\ii}[1]{\iii\{#1\}}

\begin{document}


\begin{abstract}
An explicit representation of an arbitrary zero-mean distribution as the mixture of (at-most-)two-point zero-mean distributions is given. Based in this representation, tests for (i) asymmetry patterns and (ii) for location without symmetry conditions can be constructed. Exact inequalities implying conservative properties of such tests are presented. These developments extend results established earlier by Efron, Eaton, and Pinelis under a symmetry condition. 
\end{abstract}

\title[\MakeLowercase{$t$}-test without symmetry]{{\Large Student's 
\MakeLowercase{$t$}-test without symmetry conditions}}
\date{\today; \emph{file}: \jobname.tex}
\author{Iosif Pinelis}
\address{ Department of Mathematical Sciences\\
Michigan Technological University\\
Hough\-ton, Michigan 49931 }
\email{ipinelis@mtu.edu}
\maketitle

\tableofcontents


\makeatletter
\makeatother

\section{Introduction} \label{sec:intro}


Efron \cite{efron} considered the so-called self-normalized sum
\begin{equation}\label{eq:S}
S:=\frac{X_1+\dots+X_n}{\sqrt{X_1^2+\dots+X_n^2}},
\end{equation}
assuming that the $X_i$'s are any random variables (r.v.'s) satisfying the orthant symmetry condition:
the joint distribution of $\eta_1X_1,\dots,\eta_n X_n$ is the same for any choice of signs $\eta_1,\dots,\eta_n$ in the set $\{1,-1\}$, so that, in particular, each $X_i$ is symmetric(ally distributed). It suffices that the $X_i$'s be independent and symmetrically (but not necessarily identically) distributed. 
On the event $\{X_1=\dots=X_n=0\}$, $S:=0$.


Following Efron \cite{efron}, note that the conditional distribution of any symmetric r.v.\ $X$ given $|X|$ is the symmetric distribution on the (at-most-)two-point set $\{|X|,-|X|\}$. 
Therefore, under the orthant symmetry condition, the distribution of $S$
is the mixture of the distributions of the normalized Khinchin-Rademacher sums 
$\vp_1a_1+\dots+\vp_na_n$, where the $\vp_i$'s are 
independent Rademacher r.v.'s, with $\P(\vp_i=1)=\P(\vp_i=-1)=\frac12$ for all $i$, which are also
independent of the $X_i$'s, and $a_i=X_i/(X_1^2+\dots+X_n^2)^{\frac12}$, so that $a_1^2+\dots+a_n^2=1$ (except on the event $\{X_1=\dots=X_n=0\}$, where $a_1=\dots=a_n=0$).

Let $Z\sim N(0,1)$. 
Let $a_1,\dots,a_n$ be any real numbers such that 
$a_1^2+\dots+a_n^2=1.$ The sharp form, 
\begin{equation}\label{eq:khin}
\E f\lp\vp_1a_1+\dots+\vp_n a_n\rp\le\E f(Z),
\end{equation}
of Khinchin's inequality \cite{kh} for $f(x)\equiv |x|^p$ 
was proved by
Whittle (1960) \cite{whittle} for $p\ge3$ and Haagerup (1982) \cite{haag} for $p\ge2$. 
For $f(x)\equiv e^{\la x}$ ($\la\ge0$), inequality \eqref{eq:khin} follows from Hoeffding (1963) \cite{hoeff}, whence
\begin{equation}\label{eq:exp}
\P\lp\vp_1a_1+\dots+\vp_n a_n\ge x\rp\le 
\inf_{\la\ge0}\frac{\E e^{\la Z}}{e^{\la x}}
=e^{-x^2/2}\quad \forall x\ge0.
\end{equation}

As noted by Efron, inequalities \eqref{eq:khin} and \eqref{eq:exp} together with the mentioned mixture representation imply
\begin{equation}\label{eq:khin-V}
\E e^{\la S}\le\E e^{\la Z}\quad\forall\la\ge0
\end{equation}
and
\begin{equation}\label{eq:exp-V}
\P\lp S\ge x\rp\le e^{-x^2/2}\quad \forall x\ge0.
\end{equation}
These results can be easily restated in terms of Student's statistic $T$, which is a monotonic function of $S$, as noted by Efron: $T=\sqrt{\frac{n-1}n}\,S/\sqrt{1-S^2/n}$.

Eaton (1970) \cite{eaton1} proved the Khinchin-Whittle-Haagerup inequality \eqref{eq:khin} for a much richer class of moment functions, which essentially coincides with the class $\FF3$ of all convex functions $f$ with a convex second derivative $f''$; see \cite[Proposition~A.1]{pin94} and also \cite{asymm}.
Based on this extension of \eqref{eq:khin}, inequality \eqref{eq:exp} was improved in \cite{eaton1,eaton2,pin94}. In particular,
Pinelis (1994) \cite{pin94} obtained the following improvement of a conjecture by Eaton (1974) \cite{eaton2}: 
\begin{equation*}
\P\lp\vp_1a_1+\dots+\vp_n a_n\ge x\rp\le\frac{2e^3}9\,\P(Z\ge x)\quad
\forall x\in\R. 
\end{equation*} 

Thus, inequalities \eqref{eq:khin-V} and \eqref{eq:exp-V} can be improved as follows: 
\begin{equation}\label{eq:khin-V-pin94}
\E f(S)\le\E f(Z)\quad\forall f\in\FF3
\end{equation} 
and
\begin{equation}\label{eq:exp-V-pin94}
\P\lp S\ge x\rp\le\frac{2e^3}9\,\P(Z\ge x)\quad
\forall x\in\R.
\end{equation}
Multivariate extensions of these results, which can be expressed in terms of Hotel\-ling's statistic in place of Student's, were also obtained in \cite{pin94}. 

It was pointed out in \cite[Theorem~2.8]{pin94} that, since the normal tail decreases fast, inequality \eqref{eq:exp-V-pin94} implies that relevant quantiles of $S$ may exceed the corresponding standard normal quantiles only by a relatively small amount, so that one can use \eqref{eq:exp-V-pin94} rather efficiently to test symmetry even for non-i.i.d.\ observations. 

Here we shall present extensions of inequalities \eqref{eq:khin-V-pin94} and \eqref{eq:exp-V-pin94} to the case when the $X_i$'s are not symmetric. (Asymptotics for large deviations of $S$ for i.i.d.\ $X_i$'s without moment conditions was obtained recently by Jing, Shao and Zhou \cite{jing-shao}.) 

Our basic idea is to represent any zero-mean, possibly asymmetric distribution as an appropriate mixture of two-point zero-mean distributions.  
Let us assume at first that a zero-mean r.v.\ $X$ has an everywhere continuous and strictly increasing distribution function (d.f.). Consider the truncated r.v.\ $X_{a,b}:=X\ii{a\le X\le b}$. 
\big(Here and in what follows, as usual, 
$\ii{\mathcal{A}}$ is the indicator of a given assertion $\mathcal{A}$, so that $\ii{\mathcal{A}}=1$ if $\mathcal{A}$ is true and $\ii{\mathcal{A}}=0$ if $\mathcal{A}$ is false.\big)
Then, for every fixed $a\in(-\infty,0]$, the function $b\mapsto\E X_{a,b}$ is continuous and increasing on the interval $[0,\infty)$ from $\E X_{a,0}\le0$ to 
$\E X_{a,\infty}>0$. Hence, for each $a\in(-\infty,0]$, there exists a unique value $b\in[0,\infty)$ such that $\E X_{a,b}=0$. Similarly, for each $b\in[0,\infty)$, there exists a unique value $a\in(-\infty,0]$ such that 
$\E X_{a,b}=0$. That is, one has a one-to-one correspondence between $a\in(-\infty,0]$ and $b\in[0,\infty)$ such that $\E X_{a,b}=0$. Denote by 
$\r:=\r_X$ the {\em reciprocating} function defined on $\R$ and carrying this correspondence, so that 
$$\E X\ii{\text{$X$ is between $x$ and $\r(x)$} }=0\quad\forall x\in\R;$$
the function $\r$ is decreasing on $\R$ and such that $\r(\r(x))=x$ $\forall x\in\R$; moreover, $\r(0)=0$. 
(Clearly, $\r(x)=-x$ for all real $x$ if the r.v.\ $X$ is symmetric.)
Thus, the set 
$\{\,\{x,\r(x)\}\colon x\in\R\,\}$ of (at-most-)two-point sets constitutes a partition of $\R$. Moreover, the two-point set $\{x,\r(x)\}$ is uniquely determined by the distance $|x-\r(x)|=|x|+|\r(x)|$ between the two points, as well as by the product $|x|\,|\r(x)|$. 
One can see that the conditional distribution of the zero-mean r.v.\ $X$ given $W:=|X-\r(X)|$ (or, equivalently, $Y:=|X|\,|\r(X)|$) is the uniquely determined zero-mean distribution on the two-point set $\{X,\r(X)\}$.
Thus, the distribution of the zero-mean r.v.\ $X$ with an everywhere positive density is represented as a mixture of two-point zero-mean distributions.
This mixture is given rather explicitly, provided that the distribution of r.v.\ $X$ is known. 

Thus, one has generalized versions of the self-normalized sum \eqref{eq:S}, which require -- instead of the symmetry of independent r.v.'s $X_i$ -- only that the $X_i$'s be zero-mean:
\begin{equation*} 
S_W:=\frac{X_1+\dots+X_n}{\frac12\sqrt{W_1^2+\dots+W_n^2}}\quad\text{and}\quad
S_{Y,\la}:=\frac{X_1+\dots+X_n}{(Y_1^\la+\dots+Y_n^\la)^{\frac1{2\la}}},
\end{equation*}
where $\la>0$,
$$W_i:=|X_i-\r_i(X_i)|\quad\text{and}\quad Y_i:=|X_i\,\r_i(X_i)|,$$
and the reciprocating function $\r_i:=\r_{X_i}$ is constructed as above, based on the distribution of $X_i$, for each $i$, so that the reciprocating functions $\r_i$ may be different from one another if the $X_i$'s are not identically distributed.
On the event $\{X_1=\dots=X_n=0\}$ (which is the same as either one of events $\{W_1=\dots=W_n=0\}$ and $\{Y_1=\dots=Y_n=0\}$),  $S_W:=0$ and $S_{Y,\la}:=0$. 
Note that $S_W=S_{Y,1}=S$ when the $X_i$'s are symmetric.
Logan \emph{et al} \cite{logan} and Shao \cite{shao} obtained limit theorems for the ``symmetric'' version of $S_{Y,\la}$ (with the reciprocating function $\r(x)\equiv-x$), whereas the $X_i$'s did not need to be symmetric.

These constructions can be extended to the general case of any zero-mean r.v.\ $X$, possibly with a d.f.\ which is not continuous or strictly increasing. Toward that end, one can use randomization (by means of a r.v.\ uniformly distributed in interval $(0,1)$) to deal with the atoms of the distribution of r.v.\ $X$, and generalized inverse functions to deal with the intervals on which the d.f.\ of $X$ is constant.  

Note that the reciprocating function $\r$ depends on the (usually unknown in statistics) distribution of the underlying r.v.\ $X$. However, if e.g.\ the $X_i$ constitute an i.i.d.\ sample, then the function $G$ defined by \eqref{eq:G(x)} can be estimated based on the sample, so that one can estimate the reciprocating function $\r$. Thus, replacing $X_1+\dots+X_n$ in the numerators of $S_W$ and $S_{Y,\la}$ by $X_1+\dots+X_n-n\theta$,
one obtains approximate pivots to be used to construct confidence intervals or, equivalently, tests for an unknown mean $\theta$. 
One can also use bootstrap to estimate the distributions of such pivots.

\section{Results} \label{sec:results}

Let $X$ be a zero-mean real-valued r.v.\ defined on a probability space $(\Omega,\Sigma,\P)$.
Let 
\begin{equation}\label{eq:G(x)}
	G(x):=
\begin{cases}
\E X\ii{X\in(0,x]} & \text{ if }x\in[0,\infty], \\
\E(-X)\ii{X\in[x,0)} & \text{ if }x\in[-\infty,0]. 
\end{cases}
\end{equation}

Note that $G(0)=0$; $G$ is non-decreasing and right-continuous on $[0,\infty)$; and $G$ is non-increasing and left-continuous on $(-\infty,0]$; in particular, $G$ is continuous at $0$. 
Moreover, the condition $\E X=0$ implies that
\begin{equation} \label{eq:m} 
G(\infty) = G(-\infty) =\tfrac12 \E |X|=: m<\infty. 
\end{equation}
Thus, $G(x)\in[0,m]$ for all $x\in[-\infty,\infty]$.\\

For $h\in[0,m]$, let
\begin{align}
x_+(h) & := \inf\{x\in[0,\infty]\colon G(x)\ge h\}, \label{eq:x+} \\
x_-(h)& :=\sup\{x\in[-\infty,0]\colon G(x)\ge h\}. \label{eq:x-}
\end{align}
Note that 
$x_+(h)\in[0,\infty)$ and $x_-(h)\in(-\infty,0]$ for all $h\in[0,m)$.

For $x\in\R$ and $u\in[0,1]$, define the \emph{reciprocating function} of r.v.\ $X$ by the formula  
\begin{equation} \label{eq:hat}
\r(x,u):=
\begin{cases}
x_-(H(x,u)) & \text{ if }x\in[0,\infty), \\
x_+(H(x,u)) & \text{ if }x\in(-\infty,0], \\
\end{cases}
\end{equation}
where
\begin{equation} \label{eq:H}
H(x,u):=
\begin{cases}
G(x-)+u\cdot(G(x)-G(x-)) & \text{ if }x\in[0,\infty), \\
G(x+)+u\cdot(G(x)-G(x+)) & \text{ if }x\in(-\infty,0]. 
\end{cases}
\end{equation}
Note that $H(x,u)$ depends on $u$ for a given value of $x$ only if $\P(X=x)\ne0$.

Let $U\colon\Omega\to\R$ be a r.v.\ uniformly distributed on the unit interval $[0,1]$ and independent of $X$.
For a real $x$, let
\begin{equation}\label{eq:Ux}
U_x:=
\begin{cases}
U & \text{ if } \P(X=x)\ne0,\\
1 & \text{ if } \P(X=x)=0.
\end{cases}
\end{equation}
Introduce the r.v.'s\ 
\begin{equation}\label{eq:W}
W:=|X-\r(X,U_X)|\quad\text{and}\quad
Y:=|X\,\r(X,U_X)|
\end{equation}
where the r.v.\ $U_X$ is defined in the usual manner: $U_X(\omega) := U_{X(\omega)}(\omega)$, for all $\omega\in\Omega$. 

\begin{theorem}\label{th:cond}

\begin{description}
\item[(i)]
There exist an event $\Omega_0\in\Sigma$ such that $\P(\Omega_0)=1$ and continuous functions $c\colon V_0\to(-\infty,0]$ and $d\colon V_0\to[0,\infty)$ defined on the set $V_0:=\{W(\omega)\colon\omega\in\Omega_0\}$ such that $d$ and $(-c)$ are nondecreasing on $V_0$, and on $\Omega_0$ one has
$$\{X,\r(X,U_X)\}=\{c(W),d(W)\}\quad\text{and}\quad
d(W)-c(W)=W.$$
\item[(ii)]
the conditional distribution of $X$ given $W$ coincides with that of $D_W$:
\begin{equation}\label{eq:cond}
\LL(X|W)=\LL(D_W|W),
\end{equation}
where, for every $v\in V_0$, $D_v$ is a r.v.\ such that
\begin{equation*} 
D_v=
\begin{cases}
d(v) & \text{ with probability }\frac{|c(v)|}{|c(v)|+d(v)},\\
c(v) & \text{ with probability }\frac{d(v)}{|c(v)|+d(v)}
\end{cases}
\end{equation*}
if $v\ne0$, and
$D_0\equiv0$,
so that $D_v$ takes on at most two distinct values and
$$\E D_v=0.$$
Formally, \eqref{eq:cond} is understood as follows: 
\begin{equation}\label{eq:cond-precise}
\E f(X)\ii{W\in B} = \E \vpi_f(W)\ii{W\in B}
\end{equation}
for all Borel functions $f\colon\R\to[0,\infty)$ and all Borel sets $B\subseteq[0,\infty)$, where
\begin{equation} \label{eq:vpi-f}
\vpi_f(v):= \E f(D_v)=
\begin{cases}
f(c(v)) \frac{d(v)}{|c(v)|+d(v)} + f(d(v)) \frac{|c(v)|}{|c(v)|+d(v)} & \text{ if }v\ne0,\\
f(0) & \text{ if }v=0.
\end{cases}
\end{equation}
That is, \eqref{eq:cond} means that
\begin{equation}\label{eq:cond-precise1}
\E f(X)\ii{W\in B} = \int_{\R}\P(W\in\d v)\,\E f(D_v)\ii{v\in B},
\end{equation}
where $f$ and $B$ are as in \eqref{eq:cond-precise}.

This understanding differs somewhat from the way in which the notion of the conditional distribution is usually understood. The above meaning is more convenient in the applications below, because \eqref{eq:vpi-f} can be generalized as follows. 

For all Borel functions $F\colon\R\times\R\to[0,\infty)$,
\begin{equation} \label{eq:F} 
\E F(X,W)=\int_{\R} \P(W\in\d v)\,\E F(D_v,v);
\end{equation}
in fact, one can write $\int_{[0,\infty)}$ instead of $\int_{\R}$ in \eqref{eq:cond-precise1} and \eqref{eq:F}, because $W\ge0$ a.s.

\end{description}
\end{theorem}

The following theorem is quite similar to Theorem~\ref{th:cond}.

\begin{theorem}\label{th:cond-y}

\begin{description}
\item[(i)]
There exist an event $\Omega_0\in\Sigma$ such that $\P(\Omega_0)=1$ and continuous functions $\tilde c\colon\tilde V_0\to(-\infty,0]$ and $\tilde d\colon\tilde V_0\to[0,\infty)$ defined on the set $\tilde V_0:=\{Y(\omega)\colon\omega\in\Omega_0\}$ such that $\tilde d$ and $(-\tilde c)$ are nondecreasing on $\tilde V_0$, and on $\Omega_0$ one has
$$\{X,\r(X,U_X)\}=\{\tilde c(Y),\tilde d(Y)\}\quad\text{and}\quad
-\tilde c(Y)\tilde d(Y)=Y.$$
\item[(ii)]
the conditional distribution of $X$ given $Y$ coincides with that of $\tilde D_Y$:
\begin{equation}\label{eq:cond-y}
\LL(X|Y)=\LL(\tilde D_Y|Y),
\end{equation}
where, for every $y\in\tilde V_0$, $\tilde D_y$ is a r.v.\ such that
\begin{equation*} 
\tilde D_y=
\begin{cases}
\tilde d(y) & \text{ with probability }
\frac{|\tilde c(y)|}{|\tilde c(y)|+\tilde d(y)},\\
\tilde c(y) & \text{ with probability }
\frac{\tilde d(y)}{|\tilde c(y)|+\tilde d(y)}
\end{cases}
\end{equation*}
if $y\ne0$, and
$\tilde D_0\equiv0$,
so that $\tilde D_y$ takes on at most two distinct values and
$$\E\tilde D_y=0.$$
Formally, \eqref{eq:cond-y} is understood as follows. 

For all Borel functions $F\colon\R\times\R\to[0,\infty)$,
\begin{equation} \label{eq:F-y} 
\E F(X,Y)=\int_{\R} \P(Y\in\d y)\,\E F(\tilde D_y,y);
\end{equation}
in fact, one can write $\int_{[0,\infty)}$ instead of $\int_{\R}$ in \eqref{eq:F-y}, because $Y\ge0$ a.s.

\end{description}
\end{theorem}

\begin{remark*}
It is easily seen from the proof of Theorem~\ref{th:cond-y} or, more specifically, from the proof of Lemma~\ref{lem:product}, that Theorem~\ref{th:cond-y} holds for all r.v.'s $Y$ of the more general form $\psi(|X|,|\r(X,U_X)|)$, where $\psi(u,v)$ is any expression such that (i) $\psi(0,0)=0$; (ii) $\psi(u,v)$ is nondecreasing in $u$ and in $v$ over all nonnegative $u$ and $v$; and (iii) $\psi(u,v)$ is strictly increasing in $u$ and in $v$ over all strictly positive $u$ and $v$.  
\end{remark*}

\begin{example}\label{ex:discrete}
Let $X$ have the discrete distribution
$\frac5{10}\,\de_{-1}+\frac1{10}\,\de_0+\frac3{10}\,\de_1+\frac1{10}\,\de_2$ on the finite set $\{-1,0,1,2\}$, where $\de_a$ denotes the (Dirac) probability distribution on the singleton set $\{a\}$. Then $m=\frac5{10}$ and, for $x\in\R$, $u\in(0,1)$, and $h\in[0,m]$,
\begin{gather*}
G(x)=\tfrac5{10}\ii{x\le-1}+\tfrac3{10}\ii{1\le x<2}+\tfrac5{10}\ii{2\le x},\\
x_+(h)=\ii{0<h\le\tfrac3{10}}+2\ii{\tfrac3{10}<h},\quad
x_-(h)=-\ii{0<h},\\
H(-1,u)=\tfrac5{10}\,u,\quad
H(0,u)=0,\quad
H(1,u)=\tfrac3{10}\,u,\quad
H(2,u)=\tfrac3{10}+\tfrac2{10}\,u,\\
\r(-1,u)=\ii{u\le\tfrac35}+2\ii{u>\tfrac35},\ \r(0,u)=0,\ \r(1,u)=-1,\ \r(2,u)=-1.
\end{gather*}
Therefore, the distribution of $W$ is $\frac1{10}\,\de_0+\frac6{10}\,\de_2+\frac3{10}\,\de_3$ and the conditional distributions of $X$ given $W=0$, $W=2$, and $W=3$ are $\de_0$, $\frac12\,\de_{-1}+\frac12\,\de_1$, and $\frac23\,\de_{-1}+\frac13\,\de_2$, respectively.
Thus, the zero-mean distribution of $X$ is represented as a mixture of (at-most-)two-point zero-mean distributions:
$$
\tfrac5{10}\,\de_{-1}+\tfrac1{10}\,\de_0+\tfrac3{10}\,\de_1+\tfrac1{10}\,\de_2
=
\tfrac1{10}\,\de_0
+\tfrac6{10}\,(\tfrac12\,\de_{-1}+\tfrac12\,\de_1)
+\tfrac3{10}\,(\tfrac23\,\de_{-1}+\tfrac13\,\de_2).
$$ 
Equivalently, one can condition here on $Y$ instead of $W$. 
The distribution of $Y$ is $\frac1{10}\,\de_0+\frac6{10}\,\de_1+\frac3{10}\,\de_2$ and the conditional distributions of $X$ given $Y=0$, $Y=1$, and $Y=2$ are $\de_0$, $\frac12\,\de_{-1}+\frac12\,\de_1$, and $\frac23\,\de_{-1}+\frac13\,\de_2$, respectively.
\end{example}

\begin{remark*}
A zero-mean distribution can be represented as a mixture of (at-most-)\\
two-point zero-mean distributions in a variety of ways. 
For instance, the symmetric distribution 
$\frac1{10}\de_{-2}+\frac4{10}\de_{-1}+\frac4{10}\de_{1}+\frac1{10}\de_{2}$ can be represented either as the mixture
$\frac3{10}(\frac13\de_{-2}+\frac23\de_{1})+
\frac3{10}(\frac13\de_{2}+\frac23\de_{-1})+
\frac4{10}(\frac12\de_{-1}+\frac12\de_{1})$ of two asymmetric and one symmetric two-point zero-mean distributions or as the mixture
$\frac15(\frac12\de_{-2}+\frac12\de_{2})+
\frac45(\frac12\de_{-1}+\frac12\de_{1})$ of two symmetric two-point zero-mean distributions.
The latter, ``more symmetric'' representation coincides with the one produced by the method of Theorem~\ref{th:cond} (or, equivalently, by that of Theorem~\ref{th:cond-y}). 
It appears that in general this method will produce the mixture representation that is ``the most symmetric'' in an appropriate sense, and hence the best with respect to such applications as Corollaries~\ref{cor:student-normal} and \ref{cor:stud-asymm}, given below. 
\end{remark*}

Let us now apply Theorems~\ref{th:cond} and \ref{th:cond-y} to the mentioned asymmetry-corrected versions of self-normalized sums.

\begin{theorem}\label{th:student}
Suppose that $X_1,\dots,X_n$ are independent zero-mean r.v.'s
and \\
$U_1,\dots,U_n$ are independent r.v.'s uniformly distributed on $[0,1]$, which are also independent of $X_1,\dots,X_n$.
For each $i=1,\dots,n$, let
$$W_i:=|X_i-\r_i(X_i,(U_i)_{X_i})|$$
be a r.v.\ constructed based on $X_i$ and $U_i$ the way the r.v.\ $W=|X-\r(X,U_X)|$ was constructed in \eqref{eq:Ux} and \eqref{eq:W} based on $X$ and $U$, where 
$\r_i$ is the reciprocating function for (the distribution of) r.v.\ $X_i$. 
Let
\begin{equation*} 
S_W:=\frac{X_1+\dots+X_n}{\frac12\sqrt{W_1^2+\dots+W_n^2}},
\end{equation*}
where the rule $\frac00:=0$ is used if the denominator is zero.
Then for every nonnegative Borel function $f$ on $\R$
\begin{equation} \label{eq:stud-ineq} 
\E f(S_W) \le \max\big(f(0),\,\sup\E f(Z_1+\dots+Z_n)\big),
\end{equation}
where the $\sup$ is taken over all $n$-tuples of independent zero-mean r.v.'s $Z_1,\dots,Z_n$ with the property that each $Z_i$ takes on only two values, say $c_i$ and $d_i$, such that
$$\frac12\sqrt{ (d_1-c_1)^2+\dots+(d_n-c_n)^2 }=1.$$
\end{theorem}

For every natural $\al$, let 
$\H\al$ denote
the class 
of all functions $f\colon\R\to\R$ such that $f$ has finite derivatives $f^{(0)}:=f,f^{(1)}:=f',\dots,f^{(\al-1)}$ on $\R$, $f^{(\al-1)}$ is convex on $\R$, and $f^{(j)}(-\infty+)=0$ for $j=0,1,\dots,\al-1$. 

\begin{corollary} \label{cor:student-normal}
Under the conditions of Theorem \ref{th:student},
\begin{align*} 
\E f(S_W) &\le \E f(Z)\quad\forall f\in\H5\quad\text{and} 
\\
\P(S_W\ge x)&\le c_{5,0}\P(Z\ge x)\quad\forall x\in\R, 
\end{align*}
where $c_{5,0}=5!(e/5)^5=5.699\dots$. 
\end{corollary}  

This follows immediately from Theorem \ref{th:student} and results of \cite{normal}. 
(Note that every function $f\in\H5$ is convex, and so, by Jensen's inequality, $f(0)\le\E f(Z)$.)

The following theorem is quite similar to Theorem~\ref{th:student}. 

\begin{theorem}\label{th:student-y}
With the $X_i$'s and $U_i$'s as in Theorem~\ref{th:student},
let for each $i=1,\dots,n$
$$Y_i:=|X_i\,\r_i(X_i,(U_i)_{X_i})|,$$
where $\r_i$ is the reciprocating function for r.v.\ $X_i$.
For any $\la>0$, let
\begin{equation*} 
S_{Y,\la}:=\frac{X_1+\dots+X_n}
{(Y_1^\la+\dots+Y_n^\la)^{\frac1{2\la}}},
\end{equation*}
where the rule $\frac00:=0$ is used if the denominator is zero.
Then for every nonnegative Borel function $f$ on $\R$
\begin{equation*} 
\E f(S_{Y,\la}) \le \max\big(f(0),\,\sup\E f(Z_1+\dots+Z_n)\big),
\end{equation*}
where the $\sup$ is taken over all $n$-tuples of independent zero-mean r.v.'s $Z_1,\dots,Z_n$ with the property that each $Z_i$ takes on only two values, say $c_i$ and $d_i$, such that
$$|c_1d_1|^\la+\dots+|c_nd_n|^\la=1.$$
(Note that $\Var Z_i=|c_id_i|$ for all $i$.)
\end{theorem}

\begin{corollary} \label{cor:stud-asymm}
Under the conditions of Theorem \ref{th:student-y}, 
suppose that for some $p\in(0,1)$ and all $i\in\{1,\dots,n\}$
\begin{equation}\label{eq:bounded-asymm}
\frac{X_i}{|\r_i(X_i,(U_i)_{X_i})|}\ii{X_i>0} \le \frac{1-p}p\ \text{a.s.}	
\end{equation}
Then for all 
\begin{equation}\label{eq:la_*}
\la\ge \la_*(p):=
\begin{cases}
\dfrac{1 + p + 2\,p^2}
   {2{\big( {\sqrt{p - p^2}} + 
       2\,p^2 \big) }} \quad & \text{if}\quad 0<p\le\frac{1}{2}, \\
1 \quad & \text{if}\quad \frac{1}{2}\le p<1.     
\end{cases}
\end{equation}
one has
\begin{align*} 
\E f(V_{Y,\la}) &\le \E f(T_n)\quad\forall f\in\H3\quad\text{and} 
\\
\P(V_{Y,\la}\ge x)&\le c_{3,0}\P^\lc(T_n\ge x)\quad\forall x\in\R, 
\end{align*}
where $T_n:=(Z_1+\dots+Z_n)/n^{1/(2\la)}$; 
$Z_1,\dots,Z_n$ are independent r.v.'s each having the standardized Bernoulli distribution with parameter $p$; 
the function $x\mapsto\P^\lc(T_n\ge x)$ is the least log-concave majorant of the function $x\mapsto\P(T_n\ge x)$ on $\R$; $c_{3,0}=2e^3/9=4.4634\ldots$. 
The upper bound $c_{3,0}\P^\lc(T_n\ge x)$ can be replaced by somewhat better ones, in accordance with \cite[Theorem~2.3]{binom} or \cite[(3.3)]{asymm}.
The lower bound $\la_*(p)$ on $\la$ given by \eqref{eq:la_*} is the best possible one, for each $p$.
\end{corollary} 

Condition \eqref{eq:bounded-asymm} is likely to hold when the $X_i$'s are bounded i.i.d.\ r.v.'s.
For instance, \eqref{eq:bounded-asymm} holds with $p=\frac13$ for r.v.\ $X$ in Example~\ref{ex:discrete} in place of $X_i$.

Corollary~\ref{cor:stud-asymm} follows immediately from Theorem \ref{th:student} and results of \cite{asymm}. 

\section{Proofs} \label{sec:proofs}

We shall precede the proof of the theorems by the statements of a number of lemmas (in Subsection \ref{subsec:lemmata}). 
Next, we shall prove the theorems (in Subsection \ref{subsec:th proofs}).  
Finally, we shall prove the lemmas (in Subsection \ref{subsec:lem proofs}).

\subsection{Statements of lemmata} \label{subsec:lemmata} 

\numberwithin{equation}{subsection}

Without loss of generality, one may assume that in Theorem \ref{th:cond} 
$$\P(X=0)\ne1.$$
Hence, 
$$m\in(0,\infty).$$
To state our lemmas, we need to introduce more notation. Consider the sets
\begin{align*}
M_+ & :=\{x\in(0,\infty)\colon\forall y<x\ \P(X\in(y,x])>0\}, \\ 
N_+ & :=\{x\in(0,\infty)\colon \P(X=x)=0\}, \\ 
L_+ & :=\{x\in(0,\infty)\colon\exists y<x\ \P(X\in(y,x))=0\}, \\
&=\{x\in(0,\infty)\colon\exists y\in[0,x)\ \P(X\in(y,x))=0\}, \\ 
M_- & :=\{x\in(-\infty,0)\colon\forall y>x\ \P(X\in[x,y))>0\}, \\ 
N_- & :=\{x\in(-\infty,0)\colon \P(X=x)=0\}, \\ 
L_- & :=\{x\in(-\infty,0)\colon\exists y>x\ \P(X\in(x,y))=0\} \\
&=\{x\in(-\infty,0)\colon\exists y\in(x,0]\ \P(X\in(x,y))=0\}, \\  
M & :=M_+\cup M_-, \\
N & :=N_+\cup N_-, \\
L & :=L_+\cup L_-. \\
\end{align*}
Note that
\begin{equation} \label{eq:MNL}
N_+\cap L_+=(0,\infty)\setminus M_+,\quad N_-\cap L_-=(-\infty,0)\setminus M_-.
\end{equation}

Now we can introduce the sets
\begin{subequations} \label{eq:G+}
\begin{align}
\G_+ :=\Bigl\{ (x,u)\colon & \notag \\
	& x\in M_+ ,\ 0\le u\le1, \label{eq:G+1} \\
	& x\in N_+\implies u=1, \label{eq:G+2} \\
	& x\in L_+\implies u>0, \label{eq:G+3} \\
	& \P(X>x)=0 \implies (x\notin N_+\ \&\ u<1)\ \Bigr\}, \label{eq:G+4}
\end{align} 
\end{subequations}
\begin{subequations} \label{eq:G-}
\begin{align}
\G_- :=\Bigl\{ (x,u)\colon  \notag \\
	& x\in M_- ,\ 0\le u\le1, \label{eq:G-1} \\
	& x\in N_-\implies u=1, \label{eq:G-2} \\
	& x\in L_-\implies u>0, \label{eq:G-3} \\
	& \P(X<x)=0 \implies (x\notin N_-\ \&\ u<1)\ \Bigr\}, \label{eq:G-4} 
\end{align}
\end{subequations}
\begin{equation} \label{eq:G}
\G:=\G_+\cup\G_-.
\end{equation}
Note that 
\begin{equation} \label{eq:G+-}
\G_+\cap\G_-=\emptyset,
\end{equation}
because $M_+\cap M_-\subseteq(0,\infty)\cap(-\infty,0)=\emptyset$ and in view of 
\eqref{eq:G+1} and \eqref{eq:G-1}. 

\begin{lemma} \label{lem:M}
$$\P(X\notin M\cup\{0\})=0.$$
\end{lemma}

\begin{lemma} \label{lem:x+-}
(Recall definitions \eqref{eq:m}, \eqref{eq:x+}, and \eqref{eq:x-}. )
For $h\in(0,m]$ 
\begin{align}
x_+(h) &= \min\{x\in(0,\infty]\colon G(x)\ge h\}; \label{eq:x+min} \\
x_-(h) &=\max\{x\in[-\infty,0)\colon G(x)\ge h\}; \label{eq:x-max}
\end{align}
\begin{gather}
G(y)<h\ \forall y\in[0,x_+(h)), \quad G(x_+(h))\ge h\ge G(x_+(h)-); \label{eq:x+,y} \\ 
G(y)<h\ \forall y\in(x_-(h),0], \quad G(x_-(h))\ge h\ge G(x_-(h)+). \label{eq:x-,y}  
\end{gather}
If, moreover, $h\in(0,m)$ then
$$ x_+(h)\in(0,\infty)\quad\text{and}\quad x_-(h)\in(-\infty,0).$$
\end{lemma} 

For 
$$h\in(0,m],$$
let
\begin{align}
u_+(h) := 
	\begin{cases}
\frac {h-G(x_+(h)-)} {G(x_+(h))-G(x_+(h)-)} & \text{ if } x_+(h)\notin N_+ , \\
1 & \text{ otherwise, }  
	\end{cases} 
\label{eq:u+(h)}\\
u_-(h) := 
	\begin{cases}
\frac {h-G(x_-(h)+)} {G(x_-(h))-G(x_-(h)+)} & \text{ if } x_-(h)\notin N_+ , \\
1 & \text{ otherwise. }  
	\end{cases}
\notag
\end{align}

\begin{lemma} \label{lem:bijection}
The formula
\begin{equation} \label{eq:map+} 
(0,m)\ni h \longmapsto (x_+(h),u_+(h))\in\G_+ 
\end{equation}
defines a {\em one-to-one} map of the interval $(0,m)$ {\em onto} $\G_+$, and the inverse map is given by the formula
\begin{equation} \label{eq:h+}
\G_+\ni(x,u) \longmapsto h_+(x,u):= G(x-)+u\cdot(G(x)-G(x-))\in(0,m). 
\end{equation} 
Similarly, the formula
\begin{equation} \label{eq:map-}  
(0,m)\ni h \longmapsto (x_-(h),u_-(h)) \in\G_-
\end{equation}
defines a {\em one-to-one} map of the interval $(0,m)$ {\em onto} $\G_-$, and the inverse map is given by the formula
\begin{equation} \label{eq:h-}
\G_-\ni(x,u) \longmapsto h_-(x,u):= G(x+)+u\cdot(G(x)-G(x+))\in(0,m). 
\end{equation} 
\end{lemma}

Note that 
\begin{equation*} 
H(x,u)=
\begin{cases}
h_+(x,u) & \text{ for } (x,u)\in\G_+, \\
h_-(x,u) & \text{ for } (x,u)\in\G_-, \\
\end{cases}
\end{equation*} 
where $H(x,u)$ is given by \eqref{eq:H}.

Now, using maps \eqref{eq:map+} and \eqref{eq:map-} and their inverses \eqref{eq:h+} and \eqref{eq:h-}, one can define a one-to-one map of $\G$ onto $\G$ 
\begin{equation} \label{eq:mapG}
\G\ni(x,u)\longleftrightarrow (\hat x,\hat u) \in\G
\end{equation} 
by formulas
\begin{equation} \label{eq:hat x,hat u} 
(\hat x,\hat u) := 
\begin{cases}
( x_-(h_+(x,u)), u_-(h_+(x,u)) ) & \text{ if } (x,u)\in\G_+, \\
( x_+(h_-(x,u)), u_+(h_-(x,u)) ) & \text{ if } (x,u)\in\G_-. \\
\end{cases}
\end{equation}
Thus, the one-to-one map \eqref{eq:mapG} is inverse to itself. It maps $\G_+$ onto $\G_-$ and $\G_-$ onto $\G_+$, and the latter two correpondences can be presented as follows:
\begin{align*}
\G_+\ni(x,u)\longleftrightarrow h=h_+(x,u)=h_-(\hat x,\hat u)\longleftrightarrow (\hat x,\hat u)\ni\G_-,\\
\G_-\ni(x,u)\longleftrightarrow h=h_-(x,u)=h_+(\hat x,\hat u)\longleftrightarrow (\hat x,\hat u)\ni\G_+.
\end{align*}

\begin{remark}\label{rem:hat x}
For $\hat x$ defined by \eqref{eq:hat x,hat u} 
and $\r$ defined by 
\eqref{eq:hat}, one has
$$\r(x,u)=\hat x$$
for all $(x,u)\in\G$.
\end{remark}

Let us now introduce the map
\begin{equation} \label{eq:w}
[0,m]\ni h\longmapsto w(h):=x_+(h)-x_-(h).
\end{equation}
Introduce also the set
\begin{equation} \label{eq:V}
V:=\{w(h)\colon h\in(0,m]\}.
\end{equation}

\begin{lemma} \label{lem:x+-w}
The functions $x_+$, $(-x_-)$, and $w$ are nonnegative and nondecreasing on $[0,m]$, and positive and left-continuous on $(0,m]$. 
\end{lemma} 

\begin{lemma} \label{lem:w,eps}
Assume that $ w(h_2)=w(h_1)+\vp $ for some $\vp\in[0,\infty)$ and some $h_1$ and $h_2$ in $[0,m]$. Then
$$
0\le x_+(h_2) - x_+(h_1) \le \vp , \quad
0\le x_-(h_1) - x_-(h_2) \le \vp . 
$$ 
\end{lemma} 

As an immediate corollary to Lemma \ref{lem:w,eps}, one obtains the following.

\begin{lemma} \label{lem:x+- cont of w}
If $ w(h_2)=w(h_1)$ for some $h_1$, $h_2$ in $[0,m]$, then $x_+(h_2)=x_+(h_1)$ and $x_-(h_2)=x_-(h_1)$.  
Thus, for $h\in[0,m]$, the values of $x_+(h)$ and $x_-(h)$ are uniquely determined by the value of $w(h)$. Moreover, there are nonnegative nondecreasing continuous real functions, say $-c$ and $d$, defined on $V\cup\{0\}$ (see \eqref{eq:V}) such that for all $h\in[0,m]$
\begin{equation} \label{eq:cd}
x_+(h)=d(w(h))\quad x_-(h)=c(w(h)),\quad\text{and}\quad
d(w(h))-c(w(h))=w(h).
\end{equation}
Furthermore, by Lemma \ref{lem:x+-w}, the functions $c$ and $d$ vanish only at $0$ and are Lipschitz with Lipschitz constants $\le1$. 
\end{lemma}

\begin{remark}\label{rem:w} 
Take any pair $(x,u)\in\G$.
It follows from Lemma~\ref{lem:bijection}, Lemma~\ref{lem:x+- cont of w}, Remark~\ref{rem:hat x}, \eqref{eq:w}, and \eqref{eq:V} that
$v:=|x-\r(x,u)|\in V$. Moreover,
$v=d(v)-c(v)$ and
\begin{itemize}
	\item 
	if $x>0$, then $x=d(v)$ and $\r(x,u)=c(v)$;
	\item 
	if $x<0$, then $x=c(v)$ and $\r(x,u)=d(v)$.
\end{itemize}
\end{remark}

\begin{lemma} \label{lem:product}
There is a strictly increasing function $\tau\colon V\cup\{0\}\to\R$ such that 
\begin{equation}\label{eq:product}
\text{$|x|\,|\r(x,u)|=\tau\big(|x-\r(x,u)|\big)$ for all $(x,u)\in\G$, and $\tau(0)=0$.}	
\end{equation}
\end{lemma} 

For $v\in V$, let
$$ h_v:=\sup\{h\in(0,m]\colon w(h)\le v \}.$$
By the definition \eqref{eq:V} of $V$, the set $\{h\in(0,m]\colon w(h)\le v \}$ is non-empty. Moreover, by Lemma \ref{lem:x+-w}, the function $w$ is left-continuous and nondecreasing on $(0,m]$. Therefore
\begin{equation} \label{eq:hv}
h_v=\max\{h\in(0,m]\colon w(h)\le v \}
\end{equation} 
and
\begin{equation} \label{eq:w(hv)} 
w(h_v)=v.
\end{equation} 

\begin{lemma} \label{lem:iff}
For any $v\in V$,
\begin{description}
\item[(i)] if $(x,u)\in\G_+$, then
\begin{multline*}
|x-\r(x,u)|\le v \iff 
\Big(x<d(v)\ \text{or}\ \big(x=d(v)\ \& \ h_+(x,u)\le h_v\big)\Big);
\end{multline*}
\item[(ii)] if $(x,u)\in\G_-$, then
\begin{multline*}
|x-\r(x,u)|\le v \iff 
\Big(x>c(v)\ \text{or}\ \big(x=c(v)\ \& \ h_-(x,u)\le h_v\big)\Big);
\end{multline*}
\end{description}
\end{lemma}

\begin{remark*}
It can be seen from the proof of Lemma \ref{lem:iff} (or otherwise) that the condition \Big($x<d(v)$ or $\big(x=d(v)\ \& \ h_+(x,u)\le h_v\big)$\Big) can be replaced by the seemingly simpler one: \big($x\le d(v)\ \& \ h_+(x,u)\le h_v$\big). However, the form used in the formulation of Lemma \ref{lem:iff} will be more convenient when Lemma \ref{lem:iff} is applied. A similar comment can be made concerning the corresponding condition in part (ii) of Lemma \ref{lem:iff}. 
\end{remark*}

\begin{lemma} \label{lem:0excep}
Let $X$ and $U_X$ be as in Theorem \ref{th:cond}. Then
$$ \P(X\ne0,(X,U_X)\notin\G)= 0.$$
\end{lemma}

\begin{lemma} \label{lem:0expec}
Let $X$ and $U_X$ be as in Theorem \ref{th:cond}. Then for all $v\in[0,\infty)$
$$ \E X\ii{W\le v} = 0;$$
recall the definition \eqref{eq:W} of $W$. 
\end{lemma}

\begin{lemma} \label{lem:0expecB}
Let $X$ and $U_X$ be as in Theorem \ref{th:cond}. Then Lemma \ref{lem:0expec} can be generalized as follows: for any Borel set $B\subseteq[0,\infty)$, 
$$ \E X\ii{W\in B} = 0.$$
\end{lemma}


Let us say that a Borel set $C\subset(0,\infty)$ is {\em null} if $\P(W\in C)=0$. Note that, if $B$ is a null set, then  
identity \eqref{eq:cond-precise} holds, because both sides of it are zero. 

In the case when a Borel set $C\subset(0,\infty)$ is not null, it must contain a point $v\in V$. \big(Indeed, by Remark~\ref{rem:w}, 
the range of $W$ on the event $\{(X,U_X)\in\G\}$ is contained in $V$. Also, by
Lemma \ref{lem:0excep}, the event $\{X\ne0,(X,U_X)\notin\G\}$ is of zero probability. Finally, $W\in C\subset(0,\infty)$ implies $W\ne0$ and hence $X\ne0$.\big)

In the case when a bounded Borel set $C\subset(0,\infty)$ is not null, let
\begin{align*}
d_{\max}(C) & :=\sup\{d(v)\colon v\in V \cap C\}, \\
d_{\min}(C) & :=\inf\{d(v)\colon v\in V \cap C\}, \\
c_{\max}(C) & :=\sup\{c(v)\colon v\in V \cap C\}, \\
c_{\min}(C) & :=\inf\{c(v)\colon v\in V \cap C\};
\end{align*}
note that, by Lemma \ref{lem:x+- cont of w}, the first two of these four numbers are in $[0,\infty)$, while the last two of them are in $(-\infty,0]$. 

In addition, for any Borel function $f\colon\R\to[0,\infty)$, let 
\begin{align*}
f_{r,\max}(C) & :=\sup\{f(x)\colon x\in[d_{\min}(C),d_{\max}(C)]\}, \\
f_{r,\min}(C) & :=\inf\{f(x)\colon x\in[d_{\min}(C),d_{\max}(C)]\}, \\
f_{\ell,\max}(C) & :=\sup\{f(x)\colon x\in[c_{\min}(C),c_{\max}(C)]\}, \\
f_{\ell,\min}(C) & :=\inf\{f(x)\colon x\in[c_{\min}(C),c_{\max}(C)]\}.
\end{align*}
Here, $r$ and $\ell$ stand for ``right" and ``left", respectively.

For any $\vp>0$, let us say that a bounded Borel set $C\subset(0,\infty)$ is {\em $(d,\vp)$-good} if it is not null and is such that
$$  0 < d_{\max}(C)\le e^\vp d_{\min}(C)  . $$
Similarly, let us say that a bounded Borel set $C$ is {\em $(c,\vp)$-good} if it is not null and is such that
$$ 0 < -c_{\min}(C) \le e^\vp (-c_{\max}(C))  ; $$
recall that $ c_{\min}(C) \le c_{\max}(C) \le 0$, for any $C\subset(0,\infty)$.
 
Let us say that a bounded Borel set $C\subset(0,\infty)$ is {\em $(f,\vp)$-good} if it is not null and is such that
$$ 0<f_{r,\max}(C)\le e^\vp f_{r,\min}(C) \quad\text{and}\quad 0<f_{\ell,\max}(C)\le e^\vp f_{\ell,\min}(C)   . $$

Let us say that $C$ is {\em $\vp$-good} if it is $(d,\vp)$-good, $(c,\vp)$-good, and $(f,\vp)$-good.

Let us say that a partition of a bounded Borel set $B$ is \emph{Borel} if every member of the partition is a Borel set. Let us say that such a partition 
is {\em $(d,\vp)$-good} if every member set of the partition is either null or $(d,\vp)$-good. 
Similarly defined are $(c,\vp)$-good, $(f,\vp)$-good, and $\vp$-good partitions.

\begin{lemma} \label{lem:partition}
For any bounded Borel set $B\subset(0,\infty)$, any $\vp\in(0,\infty)$, and any everywhere strictly positive and continuous function $f$, there always exists an $\vp$-good partition of $B$.
\end{lemma}

\begin{lemma} \label{lem:Ef(X)I{}}
For any Borel function $f\colon\R\to[0,\infty)$, any bounded Borel set $C\subset(0,\infty)$, and any $\vp\in(0,\infty)$, if $C$ is null or $\vp$-good, then (recall \eqref{eq:vpi-f})
\begin{equation} \label{eq:Ef(X)I{}}
\begin{split}
e^{-4\vp}\,\E\vpi_f(W)\ii{W\in C}\le
\E f(X)\ii{W\in C} \\
\le
e^{4\vp}\,\E\vpi_f(W)\ii{W\in C}.
\end{split}
\end{equation}
\end{lemma}

Let 
$$D^{(1)}_{v_1},\dots,D^{(n)}_{v_n}$$
be independent r.v.'s such that, for each $j\in\{1,\dots,n\}$, the r.v.\ $D^{(j)}_{v_j}$ is constructed based on the distribution of $X_j$ the way the r.v.\ $D_v$ was constructed in Theorem \ref{th:cond} based on the distribution of $X$.

\begin{lemma} \label{lem:F}
Let $F(x_1,v_1,\dots,x_n,v_n)$ be a nonnegative Borel function of its $2n$ real arguments. 
Let $X_1,\dots,X_n,W_1,\dots,W_n$ be as in Theorem \ref{th:student}.
Then identity \eqref{eq:F} can be generalized as follows:
\begin{equation} \label{eq:Fn} 
\begin{split}
\E F(X_1,W_1, & \dots,X_n,W_n)\\
=\int_{\R^n} & \lp\prod_{i=1}^n \P(W_i\in\d v_i)\rp\, 
\E F(D^{(1)}_{v_1},v_1,\dots,D^{(n)}_{v_n},v_n).
\end{split}
\end{equation}
\end{lemma}

\subsection{Proofs of the theorems} \label{subsec:th proofs}

\begin{proof}[Proof of Theorem \ref{th:cond}]

{\bf (i)} Let $\Omega_0:=\{X=0\}\cup\{(X,U_X)\in\G\}$. Then, by Lemma \ref{lem:0excep}, one has $\P(\Omega_0)=1$ and, by Remark~\ref{rem:w}, $V_0\subseteq V\cup\{0\}$. The rest of part (i) of Theorem \ref{th:cond} now follows by Remark~\ref{rem:w} and Lemma~\ref{lem:x+- cont of w}. 

{\bf (ii)} Here we need to prove identities \eqref{eq:cond-precise} and \eqref{eq:F}. We shall do this in a few steps.

{\em Step 1.} Here we shall prove \eqref{eq:cond-precise} assuming that (a) the function $f$ is continuous and strictly positive everywhere on $\R$ and (b) the Borel set $B$ is a bounded subset of $(0,\infty)$.

By Lemma \ref{lem:partition}, for any $\vp\in(0,\infty)$, there exists an $\vp$-good partition of $B$.
Applying Lemma \ref{lem:Ef(X)I{}} to every member set of such a partition and then summing over all the member sets, one sees that inequalities \eqref{eq:Ef(X)I{}} hold for the entire set $B$, in place of $C$. 

Since $\vp>0$ was chosen arbitrarily, this implies that \eqref{eq:cond-precise} holds whenever the function $f$ is continuous and strictly positive everywhere on $\R$ and $B$ is a bounded Borel subset of $(0,\infty)$. Thus, Step 1 of the proof of \eqref{eq:cond-precise} is now complete. 

{\em Step 2.} If $B$ is any Borel subset of $(0,\infty)$, then the sets $B_n:=B\cap(0,n]$ are bounded for all $n\in(0,\infty)$, so that, according to Step 1, \eqref{eq:cond-precise} holds with $B_n$ in place of $B$. It remains to let $n\to\infty$ to see that \eqref{eq:cond-precise} holds whenever the function $f$ is continuous and strictly positive everywhere on $\R$ and the set $B$ is any Borel subset of $(0,\infty)$. 

{\em Step 3.} By \eqref{eq:hat}, if $x\ne0$, then $-\r(x,u)$ is either 0 or of the same sign as $x$. Hence, one always has $|W|=|X-\r(X,U_X)|\ge|X|$, so that $W=0$ always implies $X=0$. Therefore and in view of \eqref{eq:vpi-f}, identity \eqref{eq:cond-precise} holds for any function $f$ provided that $B=\{0\}$. Thus (cf. Step 2), \eqref{eq:cond-precise} holds whenever the function $f$ is continuous and strictly positive everywhere on $\R$ and the set $B$ is any Borel subset of $[0,\infty)$. 

{\em Step 4.} Since the $\sigma$-algebra generated by the set of all bounded continuous strictly positive on $\R$ functions is the entire Borel $\sigma$-algebra, we conclude by a functional form of a monotone class argument that \eqref{eq:cond-precise} holds whenever $f$ is a nonnegative Borel function on $\R$ \big(and the set $B$ is any Borel subset of $[0,\infty)$.\big)

{\em Step 5.} Identity \eqref{eq:cond-precise} (or its equivalent \eqref{eq:cond-precise1}) implies that \eqref{eq:F} holds for all Borel functions $F$ of the form $F(x,v)=\ii{x\in A,v\in B}$. Then, again by a monotone class argument, \eqref{eq:F} continues to hold for all nonnegative Borel functions $F$. 

The proof of Theorem \ref{th:cond} is now complete. 
\end{proof}

\begin{proof}[Proof of Theorem \ref{th:cond-y}]
Take here the same $\Omega_0$ as in the proof of Theorem~\ref{th:cond}. Then, by Lemma~\ref{lem:product}, on $\Omega_0$ the r.v.\ $Y$ is a strictly increasing (and hence one-to-one) transformation $\tau$ of r.v. $W$. 
Now Theorem \ref{th:cond-y} follows, with 
$\tilde c:=c\circ\tau^{-1}$ and $\tilde d:=d\circ\tau^{-1}$.
\end{proof}

\begin{proof}[Proof of Theorem \ref{th:student}]
The idea of the proof is simple. Since $X_1,\dots,X_n,\\
U_1,\dots,U_n$ are all independent and, for each $i$, the r.v.\ $W_i$ is a function of $X_i$ and $U_i$, it follows that the pairs $(X_1,W_1),\dots,(X_n,W_n)$ are independent. Therefore, for each $i$, the conditional distribution of $X_i$ given $W_1,\dots,W_n$ is the same as that of $X_i$ given $W_i$. By Theorem \ref{th:cond}, the latter conditional distribution coincides a.s.\ with the unique zero-mean distribution on the set $\{c_i(W_i),d_i(W_i)\}$, where the functions $c_i$ and $d_i$ are constructed based on the (original, unconditional) distribution of $X_i$ the way the functions $c$ and $d$ were constructed in the proof of part (i) of Theorem \ref{th:cond} based on the distribution of $X$; at that, $d_i(W_i)-c_i(W_i)=W_i$ a.s. Hence, 
conditionally on $W_1,\dots,W_n$, the r.v.'s 
$$\tilde Z_i:=\frac{X_i}{\frac12\sqrt{W_1^2+\dots+W_n^2}},\quad i=1,\dots,n,$$
are independent and each $\tilde Z_i$ is zero-mean and takes on (at most) two values, $$\tilde c_i:=\frac{c_i(W_i)}{\frac12\sqrt{W_1^2+\dots+W_n^2}} \quad\text{and}\quad
 \tilde d_i:=\frac{d_i(W_i)}{\frac12\sqrt{W_1^2+\dots+W_n^2}},$$
so that $\tilde d_i-\tilde c_i=W_i/(\frac12\sqrt{W_1^2+\dots+W_n^2})$ a.s., whence a.s.  
$$\frac12\sqrt{ \sum_{i=1}^n(\tilde d_i-\tilde c_i)^2}=1.$$
This implies that, for all nonnegative Borel functions $f$ 
$$\E(f(S_W)|W_1,\dots,W_n) \le \max\big(f(0),\sup\E f(Y_1+\dots+Y_n)\big)$$
a.s., where the $\sup$ is described in the statement of Theorem \ref{th:student}. 
Now inequality \eqref{eq:stud-ineq} follows.

Let us now give a formal proof of this inequality; it is based on Lemma~\ref{lem:F}.

Since $W_i\ge0$ a.s.\ for all $i=1,\dots,n$, integral $\int_{\R^n}$ in \eqref{eq:Fn} can be replaced by $\int_{[0,\infty)^n}$. Therefore, under the conditions of Lemma \ref{lem:F}, one has the inequality
\begin{equation} \label{eq:FnSup} 
\begin{split}
\E F(X_1,W_1, & \dots,X_n,W_n)\\
\le \sup\{ & 
\E F(D^{(1)}_{v_1},v_1,\dots,D^{(n)}_{v_n},v_n)\colon
(v_1,\dots,v_n)\in[0,\infty)^n \}.
\end{split}
\end{equation}

Now, for any nonnegative Borel function $f$ on $\R$, let
$$F_f(x_1,v_1,\dots,x_n,v_n):=
\begin{cases}
f\lp\dfrac{x_1+\dots+x_n}{\frac12\sqrt{v_1^2+\dots+v_n^2}}\rp & \text{ if } v_1^2+\dots+v_n^2\ne0, \\
f(0) & \text{ otherwise. }
\end{cases}
$$
Note that, for $i=1,\dots,n$,  
the r.v.'s 
$$\tilde Z_i:=
\begin{cases}
\dfrac{D^{(n)}_{v_n}}{\frac12\sqrt{v_1^2+\dots+v_n^2}} & \text{ if } v_1^2+\dots+v_n^2\ne0, \\
0 & \text{ otherwise }  \\
\end{cases}$$
are independent, and each $\tilde Z_i$ is zero-mean and -- provided that $v_1^2+\dots+v_n^2\ne0$ -- takes on (at most) two values, $$\tilde c_i:=\frac{c_i(v_i)}{\frac12\sqrt{v_1^2+\dots+v_n^2}} \quad\text{and}\quad
 \tilde d_i:=\frac{d_i(v_i)}{\frac12\sqrt{v_1^2+\dots+v_n^2}},$$
so that 
$\tilde d_i-\tilde c_i=v_i/(\frac12\sqrt{v_1^2+\dots+v_n^2})$ 
and
$$\frac12\sqrt{ \sum_{i=1}^n(\tilde d_i-\tilde c_i)^2}=1.$$
This and inequality \eqref{eq:FnSup} imply inequality \eqref{eq:stud-ineq} for all nonnegative Borel functions $f$.  
\end{proof}

\begin{proof}[Proof of Theorem \ref{th:student-y}]
This proof is quite similar to that of Theorem \ref{th:student}, using Theorem~\ref{th:cond-y} in place of Theorem~\ref{th:cond}.
\end{proof}

\subsection{Proofs of the lemmata} \label{subsec:lem proofs}

\begin{proof}[Proof of Lemma \ref{lem:M}]
For every $x\in\R\setminus(M\cup\{0\})$, let $\Delta_x$ denote the union of the set, say $\J_x$, of all (closed, open, or semi-open) intervals $\de$ such that $\de\ni x$ and $\P(X\in\de)=0$. Then $\Delta_x$ is an interval. (Indeed, if $x_1$ and $x_2$ are in $\Delta_x$, then $x_1\in\de_1\subseteq\Delta_x$ and $x_2\in\de_2\subseteq\Delta_x$ for some intervals $\de_1\in\J_x$ and $\de_2\in\J_x$; it follows that the union $\de_1\cup\de_2$ is an interval which
is an element of the set $\J_x$, and also $\de_1\cup\de_2\supseteq\{x_1,x_2\}$.
Thus, for every two points $x_1$ and $x_2$ which are in $\Delta_x$, all the points between $x_1$ and $x_2$ are also in $\Delta_x$, so that $\Delta_x$ is an interval.) 
Moreover, the interval $\Delta_x$ is 
non-empty and, furthermore, 
it is of nonzero length, because, by the definition of $M$, for every $x\in\R\setminus(M\cup\{0\})$, the interval $\Delta_x$ contains an interval of the form $(y,x]$ for some $y<x$ or of the form $[x,y)$ for some $y>x$. 

Observe next that, for every $x\in\R\setminus(M\cup\{0\})$, one has $\P(X\in\Delta_x)=0$. Indeed, assuming that $x\in\R\setminus(M\cup\{0\})$, let $[a,b]$ be any closed subinterval of $\Delta_x$. Then there exist intervals $\de_a$ and $\de_b$ in $\J_x$ such that $a\in\de_a$ and $b\in\de_b$. Hence, $x\in\de_a\cap\de_b$, $\P(X\in\de_a)=0$, and $ \P(X\in\de_b) =0$, so that $[a,b]\subseteq\de_a\cup\de_b$, which implies $\P(X\in[a,b]) \le \P(X\in\de_a) + \P(X\in\de_b)=0$. Thus, $\P(X\in[a,b])=0$ for every closed subinterval $[a,b]$ of $\Delta_x$. If the interval $\Delta_x$ is itself closed, this implies that $\P(X\in\Delta_x)=0$. If, for instance, $\Delta_x$ is a (necessarily non-empty) interval $[c,d)$, semi-open on the right, and $d_n\uparrow d$, then
$\P(X\in\Delta_x)=\lim_n\P(X\in[c,d_n])=0$. The cases when the interval $\Delta_x$ is open or semi-open on the left are considered similarly. 
This proves the observation.

Observe further that, for any two points $x$ and $y$ in $\R\setminus(M\cup\{0\})$, the intervals $\Delta_x$  and $\Delta_y$ are either disjoint or the same. Indeed, suppose that (i) $\Delta_x$  and $\Delta_y$ are not disjoint and (ii) $\Delta_y\setminus\Delta_x \ne\emptyset$ (for instance). Then $\Delta:=\Delta_x\cup\Delta_y\in\J_x$, while $\Delta\not\subseteq\Delta_x$; this contradicts the definition of $\Delta_x$.

Therefore, the set $\{\Delta_x\colon x\in\R\setminus(M\cup\{0\})\}$ coincides (for some index set $I$) with a set $\{\de_i\colon i\in I\}$ of intervals of nonzero length such that $\de_i\cap\de_j=\emptyset$ for any two different indices $i$ and $j$ in $I$. For every $i\in I$, one can choose a rational point $r_i\in\de_i$, and these points will necessarily be distinct, since the intervals $\de_i$ are disjoint. Therefore, the index set $I$ must be countable. Since $x\in\Delta_x$ for every $x\in\R\setminus(M\cup\{0\})$,
one concludes that
\begin{multline*}
0\le \P(X\notin M\cup\{0\} ) \le \P\lp X\in\bigcup_{x\in \R\setminus M\cup\{0\} }\Delta_x\rp \\
= \P\lp X\in\bigcup_{i\in I}\de_i\rp = \sum_{i\in I}\P(X\in\de_i) =0, 
\end{multline*}
because each $\de_i$ coincides with some of the $\Delta_x$'s.
Now Lemma \ref{lem:M} follows.
\end{proof}

\begin{proof}[Proof of Lemma \ref{lem:x+-}]
Let $h\in(0,m]$. Since $m=G(\infty)=\lim_{x\uparrow\infty}G(x)$, there exists some $x\in(0,\infty]$ such that $G(x)\ge h$. For any such $x$, \eqref{eq:x+} implies $x_+(h)\le x$. Moreover, the right-continuity of $G$ on $[0,\infty)$ implies $G(x_+(h))\ge h$ (the latter inequality is trivial if $x_+(h)=\infty$). The inequality $G(x_+(h))\ge h$, together with $h>0$ and $G(0)=0$, yields $x_+(h)\ne0$. Thus, one has
\eqref{eq:x+min}, which, in turn, implies \eqref{eq:x+,y}. Relations \eqref{eq:x-max} and \eqref{eq:x-,y} are verified similarly. 
The last sentence in Lemma \ref{lem:x+-} is now obvious.
\end{proof}

\begin{proof}[Proof of Lemma \ref{lem:bijection}]

{\bf (I)} Take any $h\in(0,m)$.
At this point, let us check that $(x_+(h),u_+(h))\in\G_+$. In other words, let us check that requirements \eqref{eq:G+} are satisfied if $x$ and $u$ are replaced there by $x_+(h)$ and $u_+(h)$, respectively.

{\bf (I)(i)} Here we shall check that requirement \eqref{eq:G+1} is satisfied if $x$ and $u$ are replaced there by $x_+(h)$ and $u_+(h)$, respectively. That $0\le u_+(h)\le1$ follows immediately from \eqref{eq:u+(h)} and the second part of \eqref{eq:x+,y}. 

It remains at this point to check that $x_+(h)\in M_+$. By Lemma \ref{lem:x+-},
$x_+(h)\in(0,\infty)$. Assuming now that $x_+(h)\notin M_+$, one has $\P(X\in(y,x_+(h)])=0$ for some $y\in[0,x_+(h))$, so that $G(y)=G(x_+(h))-\E X\ii{y<X\le x_+(h)}=G(x_+(h))\ge h$, which contradicts the first part of \eqref{eq:x+,y}. Thus, requirement \eqref{eq:G+1} is checked. 

{\bf (I)(ii)} It follows immediately from \eqref{eq:u+(h)} that requirement \eqref{eq:G+2} is satisfied if $x$ and $u$ are replaced there by $x_+(h)$ and $u_+(h)$  respectively.

{\bf (I)(iii)} Here we shall check condition \eqref{eq:G+3} for $x_+(h)$ and $u_+(h)$ in place of $x$ and $u$. In view of point (I)(ii) above, one may assume that $x_+(h)\in L_+\setminus N_+$ but $u_+(h)=0$. Then $\exists y\in[0,x_+(h))$ $\P(X\in(y,x_+(h)))=0$, and
\eqref{eq:u+(h)} implies that $G(x_+(h)-)=h$. Hence, $G(y)=G(x_+(h)-)-\E X\ii{X\in(y,x_+(h))}=G(x_+(h)-)=h$, which contradicts the first part of \eqref{eq:x+,y}. 

{\bf (I)(iv)} Let us now check condition \eqref{eq:G+4} for $x_+(h)$ and $u_+(h)$ in place of $x$ and $u$. Assume that 
$\P(X>x_+(h))=0$. Then $G(x_+(h))=m$. If $x_+(h)\in N_+$, then $G(x_+(h)-)=G(x_+(h))=m>h$, which contradicts the second part of \eqref{eq:x+,y}. Hence, $x_+(h)\notin N_+$. If now $u_+(h)=1$, then \eqref{eq:u+(h)} implies $G(x_+(h))=h$, which is in a contradiction with $G(x_+(h))=m>h$. 

The verification of point (I) is now complete.

{\bf (II)} Let us check next that map \eqref{eq:map+} is onto $\G_+$. Take any $(x,u)\in\G_+$ and let
\begin{equation} \label{eq:h}
h:=G(x-)+u\cdot(G(x)-G(x-)) .
\end{equation}
We need to check that (i) $h\in(0,m)$, (ii) $x_+(h)=x$, and (iii) $u_+(h)=u$. 

{\bf (II)(i)} Here we shall check that $h\in(0,m)$. Indeed, the condition $(x,u)\in\G_+$ implies $x\in M_+$, so that 
$\P(X\in(0,x])>0$ and hence $G(x)>0$. 
If $G(x-)>0$, then \eqref{eq:h} implies $h>0$. 

Consider now the case $G(x-)=0$. Then $x\notin N_+$, because $G(x)>0$. Also, here $x\in L_+$, because the equalities 
$G(x-)=0=G(0)$ imply $\P(X\in(0,x))=0$. 
Therefore, conditions $(x,u)\in\G_+$ and \eqref{eq:G+3} imply that $u>0$, so that \eqref{eq:h} yields $h=uG(x)>0$. Thus, $h>0$ in all cases. 

It remains at this point to check that $h<m$. This follows from \eqref{eq:h} in the case $G(x)<m$, because $G(x-)\le G(x)$ and $0\le u\le1$. 
Since $G(x)\le G(\infty)=m$, it remains here to
consider the case $G(x)=m$. Then one has $\P(X>x)=0$, so that, by \eqref{eq:G+4}, $x\notin N_+$ and $u<1$. Now \eqref{eq:h} implies $h<G(x)=m$. Thus, $h<m$ in all cases. 

{\bf (II)(ii)} Here we shall check that $x_+(h)=x$. Take any $y\in[0,x)$. 
\big(Such a $y$ exists since $x\in M_+\subseteq(0,\infty)$.\big)
To obtain a contradiction, suppose that $h\le G(y)$. Then $h\le G(x-)$. On the other hand, conditions \eqref{eq:h} and $0\le u\le1$ imply $h\ge G(x-)$. Hence, $h=G(x-)$, and then \eqref{eq:h} implies $u\cdot(G(x)-G(x-))=0$, which in turn implies that either $x\in N_+$ or $x\notin L_+$ (indeed, if $x\notin N_+$, then $u\cdot(G(x)-G(x-))=0$ implies $u=0$, so that, by \eqref{eq:G+3}, one has $x\notin L_+$). Taking now \eqref{eq:MNL} into account, it follows now that $x\in M_+\cap(N_+\cup L_+^c)=L_+^c$, where we let $L_+^c:=(0,\infty)\setminus L_+$, for brevity. Hence, for all $y\in[0,x)$ such that $h\le G(y)$ one has  $\P(X\in(y,x))>0$, so that $h=G(x-)>G(y)$, a contradiction. 
Thus, $G(y)<h$ for all $y\in[0,x)$.
On the other hand, \eqref{eq:h} and $0\le u\le1$ imply $h\le G(x)$. Now \eqref{eq:x+min} yields $x_+(h)=x$. 

{\bf (II)(iii)} Here we shall check that $u_+(h)=u$. This follows from \eqref{eq:u+(h)}, \eqref{eq:h}, and (II)(ii)
in the case $x\notin N_+$. If $x\in N_+$, then, by \eqref{eq:G+2}, $u=1$, so that $u_+(h)=u$ by \eqref{eq:u+(h)}. 

The verification of point (II) is now complete.

{\bf (III)} Let us check next that map \eqref{eq:map+} is one-to-one and its inverse is given by \eqref{eq:h}. Indeed, it follows by the first line of \eqref{eq:u+(h)} in the case $x\notin N_+$ and by the second part of \eqref{eq:x+,y} in the case 
$x\in N_+$ that,
if $x_+(h)=x$ and $u_+(h)=u$, then the value of $h$ is given by \eqref{eq:h}, and is thus uniquely determined by 
$x$ and $u$.  
 
Thus, the first half of Lemma \ref{lem:bijection} is proved. The proof of its second half is quite similar. 
\end{proof}

\begin{proof}[Proof of Lemma \ref{lem:x+-w}]
That $x_+$, $-x_-$, and $w$ are nonnegative and nondecreasing on $[0,m]$ and positive on $(0,m]$ follows immediately from \eqref{eq:x+}, \eqref{eq:x-}, Lemma~\ref{lem:x+-}, and \eqref{eq:w}. 

Let now $h\in(0,m]$, $h_n\uparrow h$, $x_n:=x_+(h_n)$, and $x:=x_+(h)$. Then, because $x_+$ is nondecreasing, one has 
$x_n\nearrow y$ for some $y\in(0,x]$. 

To obtain a contradiction, assume that $y<x$. Let $z\in(y,x)$. Then, by the first part of \eqref{eq:x+,y}, $G(z)<h$. On the other hand, $y\ge x_n$ for all $n$. Hence, $h>G(z)\ge G(y)\ge G(x_n)\ge h_n$, by the second part of \eqref{eq:x+,y}. This implies 
$h>G(z)\ge h$, which is a contradiction.

It follows that $x_+$ is left-continuous on $(0,m]$; similarly, for $x_-$ and, in view of \eqref{eq:w}, for $w$. 
\end{proof}

\begin{proof}[Proof of Lemma \ref{lem:w,eps}]
To obtain a contradiction, assume that $ x_+(h_2) - x_+(h_1)\\
<0$. Then $h_2<h_1$, since $x_+$ is nondecreasing (by Lemma \ref{lem:x+-w}). Hence, again by Lemma \ref{lem:x+-w}, $ x_-(h_1) - x_-(h_2) \le0$. By re-grouping terms, it follows that
\begin{align}
0\le\vp=w(h_2)-w(h_1) & = (x_+(h_2) - x_+(h_1)) + (x_-(h_1) - x_-(h_2)) \label{eq:eps} \\
				& < 0, \notag
\end{align} 
a contradiction. 
Therefore, $ x_+(h_2) - x_+(h_1) \ge0$. Similarly is shown that $ x_-(h_1) - x_-(h_2) \ge0$. Now Lemma \ref{lem:w,eps}
follows from \eqref{eq:eps}.
\end{proof}

\begin{proof}[Proof of Lemma~\ref{lem:product}]
In view of Remark~\ref{rem:w}, the function $\tau:=|c|d$ satisfies \eqref{eq:product} (in fact, this is the only such function). By Lemma~\ref{lem:x+- cont of w}, functions $|c|$ and $d$ are nondecreasing and vanish only at $0$, and also (in view of \eqref{eq:V}) $|c|(v)+d(v)=v$ for all $v\in V\cup\{0\}\to\R$. 
It remains to show that $\tau$ is strictly increasing. 
Take any $v_1$ and $v_2$ in $V\cup\{0\}\to\R$ such that $0\le v_1<v_2$. Then $\tau(v_2)=|c|(v_2)d(v_2)>0$,
$0\le|c|(v_1)\le|c|(v_2)$, and $0\le d(v_1)\le d(v_2)$. 
So, if $|c|(v_1)=0$ or $d(v_1)=0$, then $\tau(v_1)=0<\tau(v_2)$. 
Also, the identity $|c|(v)+d(v)=v$ implies that at least one of the inequalities $|c|(v_1)\le|c|(v_2)$ and $d(v_1)\le d(v_2)$ must be strict. 
Therefore, in all cases
$\tau(v_1)=|c|(v_1)d(v_1)<|c|(v_2)d(v_2)=\tau(v_2)$.
\end{proof}

\begin{proof}[Proof of Lemma \ref{lem:iff}]
Let $v\in V$. 
Let us prove part (i) of Lemma \ref{lem:iff}. 
Accordingly, assume that $(x,u)\in\G_+$. In view of \eqref{eq:w(hv)}, \eqref{eq:cd}, and \eqref{eq:w}, one has
\begin{equation} \label{eq:x+_w(hv)}
\begin{split} 
d(v)=d(w(h_v))=x_+(h_v),\quad c(v)=c(w(h_v))=x_-(h_v), \quad \quad \ \\ 
v=w(h_v)=x_+(h_v)-x_-(h_v).
\end{split}
\end{equation}
Let now (recall \eqref{eq:h+})
\begin{equation} \label{eq:h=h+}
h:=h_+(x,u),
\end{equation}
so that, by Lemma \ref{lem:bijection} and definitions \eqref{eq:hat} and \eqref{eq:w}, 
\begin{equation} \label{eq:x+_w(h)}
\begin{split}
h\in(0,m),\quad 
x=x_+(h)>0,\quad \r(x,u)=x_-(h)<0, \quad \quad \ \\ 
|x-\r(x,u)|=x-\r(x,u)=w(h).
\end{split}
\end{equation}

Now let us prove the ``$\Longrightarrow$" implication of part (i) of Lemma \ref{lem:iff}. 
Assume that $|x-\r(x,u)|\le v$, which can be rewritten, in view of the last equality in \eqref{eq:x+_w(h)}, as $w(h)\le v$. Now it follows from \eqref{eq:hv} that 
\begin{equation} \label{eq:h,hv}
h\le h_v.
\end{equation} 
Moreover, \eqref{eq:x+_w(h)} and \eqref{eq:x+_w(hv)} together with Lemma \ref{lem:x+-w} imply that $x\le d(v)$. Thus, 
in view of \eqref{eq:h=h+} and \eqref{eq:h,hv},
the ``$\Longrightarrow$" implication is checked.
 
Next, let us prove the ``$\Longleftarrow$" implication of part (i) of Lemma \ref{lem:iff}. Indeed, consider first the case $x<d(v)$, which can be rewritten, again in view of \eqref{eq:x+_w(h)} and \eqref{eq:x+_w(hv)}, as $x_+(h)<x_+(h_v)$; then, by the ``nondecreasing" part of Lemma \ref{lem:x+-w} and \eqref{eq:w(hv)}, one has 
$h<h_v$ and hence
$|x-\r(x,u)|=w(h)\le w(h_v)=v$. 
Consider the remaining case when $x=d(v)\ \& \ h\le h_v$. Then, applying \eqref{eq:x+_w(h)}, Lemma \ref{lem:x+-w}, and \eqref{eq:w(hv)}, one obtains 
$|x-\r(x,u)|=w(h)\le w(h_v)=v$. Thus, the ``$\Longleftarrow$" implication is also checked.

Thereby, part (i) of Lemma \ref{lem:iff} is proved. Part (ii) of the lemma is proved similarly. 
\end{proof}

\begin{proof}[Proof of Lemma \ref{lem:0excep}]
Recalling the definitions of $\G$, $\G_+$, and $\G_-$ ( \eqref{eq:G}, \eqref{eq:G+}, \eqref{eq:G-}) and the relations $M_+\subseteq(0,\infty)$ and $M_-\subseteq(-\infty,0)$, one has 
\begin{equation} \label{eq:sum} 
\P(X\ne0,(X,U_X)\notin\G) = \P(X>0,(X,U_X)\notin\G_+) + \P(X<0,(X,U_X)\notin\G_-).
\end{equation}
Next,
\begin{equation} \label{eq:4sum}
\begin{split}  
0& \le\P(X>0,(X,U_X)\notin\G_+) \\
& \le\P(X>0,X\notin M_+) + \P(U=0) + \P(X\in K_+) + \P(U=1), 
\end{split}
\end{equation}
where
\begin{align*}
K_+ & :=\{x\in(0,\infty)\colon \P(X>x)=0, x\in N_+\} \\
	& =\{x\in(0,\infty)\colon \P(X\ge x)=0 \}. 
\end{align*}
The four sommands in \eqref{eq:4sum} correspond to the restrictions on $(x,u)$ in the definition of $\G_+$. Namely, the first two summands correspond to restrictions \eqref{eq:G+1} and \eqref{eq:G+3}, respectively, while the last two summands correspond to \eqref{eq:G+4}. Note that restriction \eqref{eq:G+2} is already taken care of by definition \eqref{eq:Ux} of $U_x$.

The second and the fourth summands in \eqref{eq:4sum} are zero, because r.v.\ $U$ is uniformly distributed between 0 and 1. The first summand is zero by Lemma \ref{lem:M}. If $K_+=\emptyset$, then the third summand is zero as well. 

Assume now that $K_+\ne\emptyset$. Observe that, if $x\in K_+$ and $y\in(x,\infty)$, then $0\le\P(X\ge y)\le\P(X\ge x)=0$, whence $y\in K_+$. This implies that $K_+$ is an interval, either of the form $(a,\infty)$ for some $a\in[0,\infty)$ or of the form $[a,\infty)$ for some $a\in(0,\infty)$. 

Therefore, if $a\in K_+$, then $K_+=[a,\infty)$, and so, $\P(X\in K_+)=\P(X\ge a)=0$. In the other case, when $a\notin K_+$, one has $K_+=(a,\infty)$, and so, $\P(X\in K_+)=\P(X>a)=\lim_{n\to\infty}\P(X\ge a+\frac1n)=0$. Thus, in all cases the third summand in \eqref{eq:4sum} is zero. 

Hence, $\P(X>0,(X,U_X)\notin\G_+)=0$. Similarly, $\P(X<0,(X,U_X)\notin\G_-)=0$. Now Lemma \ref{lem:0excep}
follows by \eqref{eq:sum}. 
\end{proof}

\begin{proof}[Proof of Lemma \ref{lem:0expec}]
Take any $v\in[0,\infty)$.
In view of Lemma \ref{lem:0excep} and formulas \eqref{eq:G} and \eqref{eq:G+-},
\begin{align}
\E X\, & \ii{|X-\r(X,U_X)|\le v} \notag \\
= \E & X\ii{(X,U_X)\in\G,|X-\r(X,U_X)|\le v} \label{eq:+-split}   \\
& \label{eq:2sum}
\begin{aligned}
& =\E X\ii{(X,U_X)\in\G_+,|X-\r(X,U_X)|\le v}   \\
& +\E X\ii{(X,U_X)\in\G_-,|X-\r(X,U_X)|\le v} 
\end{aligned}
\end{align}

From this point on, the proof proceeds differently depending on properties of the value of $v$. We consider separately the following cases: (I) $V\cap(0,v]=\emptyset$; (II) $v\in V$; (III) $v$ is any upper bound of $V$; and (IV) $v\in(v_1,v_2)$ for some $v_1$ and $v_2$ in $V$. These cases are clearly exhaustive. However, in general, not all of these cases are mutually exclusive.    

{\bf(I)} Consider first the case $V\cap(0,v]=\emptyset$. 
By Lemma \ref{lem:bijection} and \eqref{eq:V}, 
$(x,u)\in\G$ implies 
that $|x-\r(x,u)|\in V$. Therefore, the expression in \eqref{eq:+-split} is zero.
Thus, Lemma \ref{lem:0expec} is proved in the case $V\cap(0,v]=\emptyset$. 

{\bf(II)} Next, consider the case $v\in V$.

In this case, by Lemma \ref{lem:iff} and also again Lemma \ref{lem:0excep}, 
\begin{align}
& \E X\ii{(X,U_X)\in\G_+,|X-\r(X,U_X)|\le v}  \notag  \\
& =\E X\ii{(X,U_X)\in\G_+, \bigl(X<d(v)\text{ or }(X=d(v)\ \&\ h_+(X,U_X)\le h_v) \bigr)} \notag \\
& =\E X\ii{ 0<X<d(v)\text{ or }(X=d(v)\ \&\ h_+(d(v),U_{d(v)})\le h_v) \bigr)} \notag \\
& = G(d(v)-) + d(v) \P(X=d(v)) \P( h_+(d(v),U)\le h_v); \label{eq:!}
\end{align}
the last equality is obvious if $d(v)\in N_+$, and it follows from the definition \eqref{eq:Ux} and the independence of $X$ and $U$ if $d(v)\notin N_+$. 

Note that 
$$d(v) \P(X=d(v)) = G(d(v)) -G(d(v)-).$$
Recall that, in view of \eqref{eq:w(hv)} and \eqref{eq:cd},
\begin{equation} \label{eq:d(v)=x+}
d(v)=d(w(h_v))=x_+(h_v).
\end{equation}

If $h_v=m$ then, by Lemma~\ref{lem:x+-w} and \eqref{eq:w(hv)}, $w(h)\le w(h_v)=v$ for all $h\in[0,m]$; that is, $v$ is an upper bound of the set $V$, so that one has Case (III), to be considered next.

It remains here to consider the case $h_v<m$.

Consider the two possible subcases. 

{\em Subcase 1:} $d(v)\in N_+$. 
In view of \eqref{eq:d(v)=x+} and \eqref{eq:h+}, for any $u\in[0,1]$, the expression in \eqref{eq:!} equals
\begin{multline*}
G(d(v)-) = G(d(v)-) + u \cdot (G(d(v)) -G(d(v)-)) = h_+(d(v),u) 
= h_+( x_+(h_v) ,u). 
\end{multline*}
Now, substituting here $u_+(h_v)$ for $u$, one has by Lemma~\ref{lem:bijection} that, in Subcase 1, the expression in \eqref{eq:!} equals $h_v$. 

{\em Subcase 2:} $d(v)\notin N_+$. Here, in view of \eqref{eq:d(v)=x+}, 
\eqref{eq:h+}, and \eqref{eq:u+(h)}, one has
\begin{multline*}
\P(h_+(d(v),U)\le h_v) = \P(h_+(x_+(h_v),U)\le h_v) \\
= \P(G( x_+(h_v) -)+U\cdot (G( x_+(h_v) ) -G( x_+(h_v) -))   \le h_v) \\
= \P(U\le u_+(h_v) ) = u_+(h_v) .
\end{multline*}
Hence, in Subcase 2, the expression in \eqref{eq:!} equals 
\begin{multline*}
G(d(v)-) + u_+(h_v) \cdot (G(d(v)) -G(d(v)-)) = h_+(d(v),u_+(h_v)) \\
= h_+(x_+(h_v),u_+(h_v)) = h_v, 
\end{multline*}
by \eqref{eq:h+}, \eqref{eq:d(v)=x+}, and Lemma \ref{lem:bijection}. 

Thus, in both Subcase 1 and Subcase 2, the expression in \eqref{eq:!} equals $h_v$. That is, the first summand in \eqref{eq:2sum} equals 
$h_v$. Similarly, the second summand in \eqref{eq:2sum} equals $-h_v$. Now Lemma \ref{lem:0expec} follows -- for all $v\in V$.

{\bf(III)} Next, if $v$ is any upper bound of $V$ then, by \eqref{eq:V} and Lemma~\ref{lem:bijection},  
$(x,u)\in\G_+$ implies $|x-\r(x,u)|\le v$, so that, in view of Lemma \ref{lem:0excep}, the first summand in 
\eqref{eq:2sum} equals $\E X\ii{X>0}=m$; similarly, the second summand in 
\eqref{eq:2sum} equals $-m$. Thus, Lemma \ref{lem:0expec} is proved in the case when $v$ is any upper bound of $V$. 

{\bf(IV)} It remains to consider the case when $v\in(v_1,v_2)$ for some $v_1$ and $v_2$ in $V$, so that $v_i=w(h_i)$ for some $h_i\in(0,m]$, where $i=1,2$. Let 
\begin{equation*} 
v_{*}:=\sup (V\cap(0,v]) .
\end{equation*}
Then $v_{*} \in (0,v]$ (because $v_1\in V\cap(0,v]$ and hence $V\cap(0,v]\ne\emptyset$).

Moreover, $v_{*} \in V$, 
so that
\begin{equation} \label{eq:v**}
v_{*}:=\max (V\cap(0,v]) .
\end{equation}
Indeed, otherwise there is a strictly increasing sequence $(v_n)$ in $V\cap(0,v]$ which converges to $v_{*}$. Then, by \eqref{eq:V}, there exists a sequence $(h^*_n)$ in $(0,m]$ such that $v_n=w(h^*_n)$ for all $n$. By Lemma \ref{lem:x+-w}, the function $w$ is nondecreasing, and so, the sequence $(h^*_n)$ is necessarily increasing. 
Hence, $h:=\lim_n h^*_n
\in(0,m]$. 
Again by Lemma \ref{lem:x+-w}, the function $w$ is left-continuous on $(0,m]$, and so, $w(h)=\lim_n w(h^*_n)=\lim_n v_n=v_{*}$. Thus, the claim that $v_{*} \in V$ is checked. 

In view of Lemma \ref{lem:bijection} and \eqref{eq:V}, $(x,u)\in\G$ implies that $|x-\r(x,u)|\in V$, whence, by \eqref{eq:v**},
for all $(x,u)\in\G$,
$$ |x-\r(x,u)|\le v \iff |x-\r(x,u)| \le v_{*}.$$
Therefore and by virtue of \eqref{eq:+-split},
the case when $v\in(v_1,v_2)$ for some $v_1$ and $v_2$ in $V$ is reduced to case (II) $v\in V$.  
\end{proof}

\begin{proof}[Proof of Lemma \ref{lem:0expecB}]
Lemma \ref{lem:0expec} implies
$$ \E X\, \ii{W\in(v_1,v_2]} =  \E X\, \ii{W\le v_2} - \E X\, \ii{W\le v_1} = 0$$
for any left-open interval $(v_1,v_2]\subseteq(0,\infty)$. 
Thus, the countably additive function (c.a.f.) $B\mapsto \E X\ii{W\in B}$ is zero on the semiring of such intervals. Since this semiring generates the entire Borel $\si$-algebra in $(0,\infty)$, this c.a.f.\ is zero on this $\si$-algebra. It remains to note that $\E X\ii{W\in\{0\}} = \E X\ii{W\le0} = 0$, because $W\ge0$ a.s.\ and by Lemma \ref{lem:0expec}.
\end{proof}

\begin{proof}[Proof of Lemma \ref{lem:partition}]
For any set $A\subseteq\R$, consider its pre-images under $c$ and $d$: 
$$d^{-1}(A):=\{v\in V\colon d(v)\in A\}\quad\text{and}\quad
c^{-1}(A):=\{v\in V\colon c(v)\in A\}.
$$
Then, for any $\de\in(0,\vp]$, the sets 
$$ C_{j,k}:=
d^{-1}(( e^{\de k},e^{\de(k+1)} ]) \cap
c^{-1}([ -e^{\de(j+1)},-e^{\de j})) \cap B, $$
where $j$ and $k$ run over all integers, form a partition of $B$ which is both $(d,\vp)$-good and $(c,\vp)$-good (because, by Lemma~\ref{lem:x+- cont of w}, functions $d$ and $-c$ are (strictly) positive on $V$. 
 
It suffices to prove that this partition is also $(f,\vp)$-good, provided that $\de\in(0,\vp]$ is small enough. Toward that end, consider any one of the $C_{j,k}$'s which are not null, so that
\begin{equation} \label{eq:dmaxmin}
0<d_{\max}(C_{j,k})\le e^\de d_{\min}(C_{j,k}),
\end{equation}
by the construction of $C_{j,k}$. 

Let 
\begin{equation} \label{eq:de1}
\de_1:=(e^\vp-1)\inf\{f(x)\colon|x|\le\sup B\}.
\end{equation}
Then $\de_1>0$, 
because the set $B$ is assumed to be bounded and the function $f$, everywhere continuous and strictly positive. Then $f$ is uniformly continuous on all bounded sets, so that there exists some $\de_2>0$ such that 
$$(|x-y|\le\de_2\ \& \ |x|\le\sup B) \implies |f(x)-f(y)|\le\de_1. $$

Choose now $\de\in(0,\vp]$ to be small enough so that 
$$(e^\de-1)\sup B\le\de_2.$$

Note that
$$0\le d_{\min}(C_{j,k}) \le d_{\max}(C_{j,k}) \le d_{\max}(B) \le \sup B,$$
by Lemma \ref{lem:x+- cont of w}. 
Therefore and in view of \eqref{eq:dmaxmin},
$$0\le d_{\max}(C_{j,k})-d_{\min}(C_{j,k}) 
\le (e^\de-1)d_{\min}(C_{j,k}) \le (e^\de-1)\sup B\le\de_2,$$
and so,
$$0<f_{r,\max}(C_{j,k})\le f_{r,\min}(C_{j,k})+\de_1
\le e^\vp f_{r,\min}(C_{j,k}),$$
by the definition \eqref{eq:de1} of $\de_1$. 
Similarly,
$$0<f_{\ell,\max}(C_{j,k})\le e^\vp f_{\ell,\min}(C_{j,k}).$$ 
\end{proof}

\begin{proof}[Proof of Lemma \ref{lem:Ef(X)I{}}]
The case when $C$ is a null set is trivial, because then each of the three terms in \eqref{eq:Ef(X)I{}} is zero.

Assume now that the set $C$ is $\vp$-good.
If $(x,u)\in\G_+$, then, by Lemma~\ref{lem:bijection}, $x=x_+(h)$ for some $h\in(0,m)$; hence, by \eqref{eq:cd}, $x=d(v)$ for $v:=w(h)=x-\r(x,u)$. Therefore, if event $\{(X,U_X)\in\G_+,W\in C\}$ occurs, then $X=d(W)$, whence $X\in[d_{\min}(C),d_{\max}(C)]$.
Similarly, if event $\{(X,U_X)\in\G_-,W\in C\}$ occurs, then $X=c(W)$, whence $X\in[c_{\min}(C),c_{\max}(C)]$.
Also, if event $\{W\in C\}$ occurs, then $X\ne0$, because $C\subset(0,\infty)$, and $X=0$ implies $W=0$. 

It follows that
\begin{equation} \label{eq:inclus}
\{ W\in C \} \subseteq \{ X\in[d_{\min},d_{\max}]\cup[c_{\min},c_{\max}]\} \cup E,
\end{equation}
where
\begin{equation*} 
E:=\{X\ne0,(X,U_X)\notin\G\},
\end{equation*}
and we set, for brevity:
$$d_{\min}:=d_{\min}(C),\ d_{\max}:=d_{\max}(C),\ 
c_{\min}:=c_{\min}(C),\ c_{\max}:=c_{\max}(C).$$ 
Note that, by Lemma \ref{lem:0excep},
\begin{equation} \label{eq:E0}
\P(E)=0.
\end{equation}

It follows from Lemma \ref{lem:0expecB}, \eqref{eq:inclus}, and \eqref{eq:E0} that 
\begin{equation} \label{eq:EX le}
0=\E X\, \ii{W\in C} \le c_{\max} q + d_{\max} p,
\end{equation}
where
$$p:=\P(X\in[d_{\min},d_{\max}],W\in C),\quad q:=\P(X\in[c_{\min},c_{\max}],W\in C), $$
so that
\begin{equation} \label{eq:p+q}
p+q=\P(W\in C).
\end{equation}
Similarly, 
\begin{equation} \label{eq:EX ge}
0=\E X\, \ii{W\in C} \ge c_{\min} q + d_{\min} p,
\end{equation}
It follows from \eqref{eq:EX le} and \eqref{eq:EX ge} that 
\begin{equation} \label{eq:p,q}
\frac{-c_{\max}}{d_{\max}} q \le p \le \frac{-c_{\min}}{d_{\min}} q.
\end{equation}
Next, letting
$$f_{r,\min}:=f_{r,\min}(C),\ f_{r,\max}:=f_{r,\max}(C),\ 
f_{\ell,\min}:=f_{\ell,\min}(C),\ f_{\ell,\max}:=f_{\ell,\max}(C),$$
one has
\begin{align} 
\E f(X)\, \ii{W\in C} 
& \le f_{\ell,\max}\, q + f_{r,\max}\, p \label{eq:EfX le 1} \\
				 & \le \frac q{d_{\min}} 
				 ( f_{\ell,\max}\cdot d_{\min} + f_{r,\max}\cdot (-c_{\min}) ); \label{eq:EfX le 2} 
\end{align}
here, inequality \eqref{eq:EfX le 1} is similar to \eqref{eq:EX le}, and \eqref{eq:EfX le 2} follows from the second inequality in \eqref{eq:p,q}.

On the other hand, recalling \eqref{eq:vpi-f} and \eqref{eq:p+q}, and then also using the first inequality in \eqref{eq:p,q}, one obtains 
\begin{align} 
\E \vpi_f(W)\, \ii{W\in C} 
& \ge \frac{ f_{\ell,\min}\cdot d_{\min} + f_{r,\min}\cdot (-c_{\max}) }{ -c_{\min}+d_{\max} } \, (p+q) \notag \\
& \ge \frac{ f_{\ell,\min}\cdot d_{\min} + f_{r,\min}\cdot (-c_{\max}) } { -c_{\min}+d_{\max} } \, 
\frac q{d_{\max}}(-c_{\max}+d_{\max}) \notag \\
& = r_1 r_2 r_3 \frac q{d_{\min}} 
( f_{\ell,\max}\cdot d_{\min} + f_{r,\max}\cdot (-c_{\min}) ), 
\label{eq:vpi_f ge} 
\end{align}
where 
\begin{align*} 
r_1 & := \frac{d_{\min}}{d_{\max}} \ge e^{-\vp}, \\
r_2 & := \frac{ -c_{\max}+d_{\max} }{ -c_{\min}+d_{\max} } \ge e^{-\vp},  \\
r_3 & := \frac{ f_{\ell,\min}\cdot d_{\min} + f_{r,\min}\cdot (-c_{\max}) }{ f_{\ell,\max}\cdot d_{\min} + f_{r,\max}\cdot (-c_{\min}) } 
	\ge e^{-\vp} \frac{ d_{\min} + (-c_{\max}) }{ d_{\min} + (-c_{\min}) } 
		\ge e^{-2\vp}.   
\end{align*} 
Now it follows from \eqref{eq:vpi_f ge} and \eqref{eq:EfX le 2} that 
$$ \E \vpi_f(W)\, \ii{W\in C} \ge e^{-4\vp} \E f(X)\, \ii{W\in C} ,$$
which proves the second inequality in \eqref{eq:Ef(X)I{}}.
The first inequality in \eqref{eq:Ef(X)I{}} is proved quite similarly.
\end{proof}

\begin{proof}[Proof of Lemma \ref{lem:F}]
For $j=1,\dots,n+1$, introduce
\begin{equation} \label{eq:Phi} 
\begin{split}
\Phi_j(x_1,v_1, & \dots,x_{j-1},v_{j-1};v_j,v_{j+1},\dots,v_n) \\
&:=\E F(x_1,v_1,\dots,x_{j-1},v_{j-1},D^{(j)}_{v_j},v_j,\dots,D^{(n)}_{v_n},v_n)
\end{split}
\end{equation}
and
\begin{equation} \label{eq:Ij} 
\begin{split}
\I_j
:= \int_{\R^{n-1+j}} & \lp\prod_{i=1}^{j-1} \P(X_i\in\d x_i,W_i\in\d v_i)\rp\,
\lp\prod_{i=j}^n \P(W_i\in\d v_i)\rp\, \\
\ \quad\times
& \Phi_j(x_1,v_1,\dots,x_{j-1},v_{j-1};v_j,v_{j+1},\dots,v_n).
\end{split}
\end{equation}
Then, for all $j=1,\dots,n$,
\begin{equation} \label{eq:Ij+1} 
\I_{j+1}
= \int_{\R^{n-2+j}}  \lp\prod_{i=1}^{j-1} \P(X_i\in\d x_i,W_i\in\d v_i)\rp\,
\lp\prod_{i=j+1}^n \P(W_i\in\d v_i)\rp\, \EE_j,
\end{equation}
where
\begin{align*} 
\EE_j
&:= \E\Phi_{j+1}(x_1,v_1,\dots,x_{j-1},v_{j-1},X_j,W_j;v_{j+1},\dots,v_n) \\
&=\int_{\R} \P(W_j\in\d v_j) \E\Phi_{j+1}(x_1,v_1,\dots,x_{j-1},v_{j-1},D^{(j)}_{v_j},v_j;v_{j+1},\dots,v_n)
\\
&=\int_{\R}\, \P(W_j\in\d v_j)\,\int_{\R} \P(D^{(j)}_{v_j}\in\d x_j)\, \\ 
&\ \quad\times\Phi_{j+1}(x_1,v_1,\dots,x_{j-1},v_{j-1},x_j,v_j;v_{j+1},\dots,v_n)
\\
&=\int_{\R}\, \P(W_j\in\d v_j)\,\int_{\R} \P(D^{(j)}_{v_j}\in\d x_j)\, \\ 
&\ \quad\times\E F(x_1,v_1,\dots,x_j,v_j,D^{(j+1)}_{v_{j+1}},v_{j+1},\dots,D^{(n)}_{v_n},v_n)
\\
&=\int_{\R}\, \P(W_j\in\d v_j)\,
\E F(x_1,v_1,\dots,x_{j-1},v_{j-1},D^{(j)}_{v_j},v_j,\dots,D^{(n)}_{v_n},v_n)
\\
&=\int_{\R}\, \P(W_j\in\d v_j)\,
\Phi_j(x_1,v_1,\dots,x_{j-1},v_{j-1};v_j,\dots,v_n);
\end{align*}
the second of these 6 equalities follows by \eqref{eq:F}, and the fourth and sixth ones by \eqref{eq:Phi}.  

Now \eqref{eq:Ij+1} and \eqref{eq:Ij} imply that 
\begin{align*} 
\I_{j+1} & =
\int_{\R^{n-1+j}} \lp\prod_{i=1}^{j-1} \P(X_i\in\d x_i,W_i\in\d v_i)\rp\,
\lp\prod_{i=j}^n \P(W_i\in\d v_i)\rp\, \\
& \ \quad\quad\quad\quad\quad\times 
\Phi_j(x_1,v_1,\dots,x_{j-1},v_{j-1};v_j,v_{j+1},\dots,v_n) \\
& =\I_j,
\end{align*}
for all $j=1,\dots,n$. 
This finally implies $\I_{n+1}=\I_1$, so that
\begin{equation*} 
\begin{split}
\E F(X_1,W_1, & \dots,X_n,W_n)=\I_{n+1}\\
=\I_1=\int_{\R^n} & \lp\prod_{i=1}^n \P(W_i\in\d v_i)\rp\, 
\E F(D^{(1)}_{v_1},v_1,\dots,D^{(n)}_{v_n},v_n).
\end{split}
\end{equation*}
\end{proof}

\end{document}